\documentclass{article}

\usepackage{amsmath}
\usepackage{amsfonts}
\usepackage{amssymb}
\usepackage{amsmath}
\usepackage{amsthm}
\usepackage{latexsym}
\usepackage{graphicx,cite}

\newtheorem{lemma1}     {Lemma}[section]
\newtheorem{teorema1}   {Theorem}
\newtheorem{prop1}      [lemma1]{Proposition}
\newtheorem{coroll1}    [lemma1]{Corollary}
\newtheorem{cong1}      [lemma1]{Conjecture}
\newtheorem{remark1}    [lemma1]{Remark}
\newtheorem{defin1}     [lemma1]{Definition} 
\newtheorem{definition} [lemma1]{Definition}

\newenvironment{Lemma}[1][]
        {\begin{lemma1}[#1]\begin{samepage}}{\end{samepage}\end{lemma1}}
\newenvironment{Theorem}[1][]
        {\begin{teorema1}[#1]\begin{samepage}}{\end{samepage}\end{teorema1}}
\newenvironment{Proposition}[1][]
        {\begin{prop1}[#1]\begin{samepage}}{\end{samepage}\end{prop1}}
\newenvironment{Corollary}[1][]
        {\begin{coroll1}[#1]\begin{samepage}}{\end{samepage}\end{coroll1}}

\newenvironment{Remark}[1][]
        {\begin{remark1}[#1]\begin{samepage}}{\end{samepage}\end{remark1}}

\newcommand{\nada}[1]   {}

\newcommand{\ep}       {\varepsilon}
\newcommand{\Om}        {\Omega}

\newcommand{\bw}{\overline w}
\newcommand{\uw}{\underline w}
\newcommand{\bk}{\bar k}
\newcommand{\uk}{\underline k}
\begin{document}

\title{\textbf{Multiplicity of supercritical fronts for
    reaction-diffusion equations in cylinders}}

\author{P. V. Gordon\footnote{Department of Mathematical Sciences, New
    Jersey Institute of Technology, Newark, NJ 07102, USA} \and
  C. B. Muratov\footnotemark[\value{footnote}] \and
  M. Novaga\footnote{Dipartimento di Matematica, Universit\`a di Padova,
    Via Trieste 63, Padova, 35121, Italy}}

\date{\today}

\maketitle

\begin{abstract}
\noindent 
We study multiplicity of the supercritical traveling front solutions
for scalar reaction-diffusion equations in infinite cylinders which
invade a linearly unstable equilibrium. These equations are known to
possess traveling wave solutions connecting an unstable equilibrium to
the closest stable equilibrium for all speeds exceeding a critical
value. We show that these are, in fact, the only traveling front
solutions in the considered problems for sufficiently large speeds.
In addition, we show that other traveling fronts connecting to the
unstable equilibrium may exist in a certain range of the wave
speed. These results are obtained with the help of a variational
characterization of such solutions.
\end{abstract}

\numberwithin{equation}{section}


\section{Introduction}

Front propagation is a ubiquitous feature of many reaction-diffusion
systems in physics, chemistry and biology
\cite{murray,kapral,zeldovich,keener,fife79,merzhanov99}. A question
of particular interest that arises in a number of applications is
related to fronts that invade an unstable equilibrium (for a review,
see \cite{vansaarloos03}). One of the well-known examples of such a
system is given by a scalar reaction-diffusion equation in a
cylindrical domain:
\begin{equation}\label{pde}
u_t = \Delta u + f(u, y).
\end{equation}
Here $u = u(x, t) \in \mathbb R$ is the dependent variable, $x \in
\Sigma \subset \mathbb R^n$ is a point in a cylindrical domain $\Sigma
= \Omega \times \mathbb R$, whose cross-section $\Omega \subset
\mathbb R^{n-1}$ is a bounded domain with the boundary of class $C^2$
(not necessarily simply-connected), and $t \in \mathbb [0, +\infty)$
is time.  We will often write $x = (y, z) \in \Sigma$, where $y \in
\Omega$ is the coordinate on the cylinder cross-section and $z \in
\mathbb R$ is the coordinate along the cylinder axis. The function $f:
\mathbb R \times \Omega \to \mathbb R$ is a nonlinear reaction
term. On different connected portions of the boundary we assume either
Dirichlet or Neumann boundary conditions (see \cite{mn:cvar08} for the
motivation of this particular choice of boundary conditions)
\begin{eqnarray}
  \label{bc}
  u \bigl|_{\partial \Sigma_\pm} = 0, \qquad \nu \cdot
  \nabla u \bigl|_{\partial \Sigma_0} = 0,
\end{eqnarray}
where $\partial \Sigma_\pm = \partial \Omega_\pm \times \mathbb R$ and
$\partial \Sigma_0 = \partial \Omega_0 \times \mathbb R$. Furthermore,
we assume that $u = 0$ is an {\em unstable} solution of (\ref{pde}) in
the sense of linear stability.

By front solutions of \eqref{pde}, one understands traveling wave
solutions, which are special solutions of \eqref{pde} of the form
$u(y, z, t) = \bar u(y, z - ct)$, connecting monotonically two
distinct equilibria, i.e., stationary $z$-independent solutions of
\eqref{pde}. In the present context, of special interest are the
fronts which invade the $u = 0$ equilibrium. Assuming the invasion
happens from left to right, the profile $\bar u > 0$ of such fronts is
described by the following equation:
\begin{eqnarray}\label{trav}
  \Delta \bar u + c \bar u_z + f(\bar u, y) = 0, \qquad \lim_{z
    \to +\infty} \bar u(\cdot, z) = 0~~\mathrm{uniformly~in~} \Omega, 
\end{eqnarray}
with the same boundary conditions as in (\ref{bc}), where $c > 0$ is
the propagation speed. These solutions are known to play an important
r\^ole in the long-time behavior of solutions of \eqref{pde} (for an
overview, see e.g., \cite{berestycki02cpam}).

The problem of existence and uniqueness of solutions to \eqref{trav}
has a long history, going back to the pioneering works of Kolmogorov,
Petrovsky and Piskunov \cite{kolmogorov37}, and within the present
context has been extensively studied by Beretsycki and Nirenberg
\cite{berestycki92}, and Vega \cite{vega93cpde,vega93}. Specifically,
Berestycki and Nirenberg established existence, monotonicity, and
uniqueness of solutions of \eqref{trav} connecting zero to the
smallest positive equilibrium for all $c \geq c_1^*$, for some $c_1^*
> 0$ (for the precise definition, see Sec. \ref{s:mr}), in the case of
Neumann boundary conditions \cite[Theorem 1.9]{berestycki92}. Their
methods were later extended by Vega to problems with Dirichlet
boundary conditions \cite[Theorem 3.4(B)]{vega93cpde}. We note that in
these works uniqueness of the front solutions is understood in the
class of functions connecting zero ahead of the front to a {\em
  prescribed} equilibrium behind the front.  In general, front
solutions of \eqref{trav} may fail to be unique, even for a fixed
value of $c$, if they are allowed to approach different
equilibria. One example of multiplicity of the front solutions to
\eqref{trav} was given by Vega under specific assumptions on the set
of equilibria of \eqref{pde} \cite[Remark 3.5(B)]{vega93cpde}.

The question of multiplicity of solutions of \eqref{trav} and
\eqref{bc} without prescribing the asymptotic behavior behind the
front is largely open. The main goal of this paper is to provide some
easily verifiable conditions to establish existence or non-existence
of multiple solutions of \eqref{trav}. To this end, we develop
variational tools to study existence and multiplicity of supercritical
fronts, i.e., the solutions of \eqref{trav} with speed $c > c^*$,
where $c^* > 0$ is some critical propagation velocity, whose precise
value will be specified in Sec. \ref{s:mr}. One of the by-products of
our methods is an independent proof of existence of solutions to
\eqref{trav} obtained in \cite{berestycki92,vega93cpde}. However, it
is worth noticing that our method can yield solutions of \eqref{trav},
\eqref{bc} that are {\em different} from those of
\cite{berestycki92,vega93cpde} under easily verifiable conditions,
thus establishing multiplicity of supercritical fronts for certain
ranges of propagation speed. At the same time, our variational
characterization, together with some refined analysis of the decay of
solutions ahead of the front also yields global uniqueness of
solutions of \eqref{trav} and \eqref{bc} for all $c > c^\sharp$, for
some $c^\sharp \geq c^*$ which is readily computable. Thus, the
solutions obtained in \cite{berestycki92,vega93cpde} are the only
positive traveling wave solutions for \eqref{pde}, \eqref{bc} invading
the unstable equilibrium $u=0$ for $c$ large enough, irrespectively of
the asymptotic behavior at $z=-\infty$.  This is one of the main
contributions of this paper, and we illustrate it by a simple example
in which all the estimates are explicit.

Our approach is mostly variational, based on the ideas developed in
\cite{heinze89,heinze,m:dcdsb04,lmn:arma08,mn:cms08}, relying on the
fact that \eqref{trav} is the Euler-Lagrange equation for a suitable
functional $\Phi_c$ defined on the exponentially weighted Sobolev
space $H^1_c(\Sigma)$ (see Section \ref{s:mr} for precise
definitions).  However, in the application of the direct method of
calculus of variations to get existence of the traveling wave
solutions of interest, one runs into the difficulty that these
solutions do not belong to $H^1_c(\Sigma)$.  To overcome this
difficulty, we first construct, by non-variational methods, an
auxiliary function which solves (\ref{trav}) for $z$ sufficiently
large and subtract its contribution from the integrand of $\Phi_c$,
thus removing the divergence of the integral.  We then minimize the
modified functional to obtain a non-trivial $H_c^1(\Sigma)$ correction
to the prescribed auxiliary function.  As a direct consequence, the
existence of a traveling front is established.

We notice that the solution constructed in this way may depend on the
choice of the auxiliary function and thus may be non-unique.
Nevertheless, we show that this does not happen for sufficiently large
values of $c$, in which case the obtained solution is the unique front
solution of \eqref{trav}.  The proof of uniqueness also has both a
variational and a non-variational parts: first, using delicate
exponentially weighted $L^2$-estimates, we show that the difference of
any two solutions with the same asymptotic decay at $z=+\infty$
belongs to $H^1_c(\Sigma)$, and then we prove uniqueness by convexity
of the respective functional for $c$ large enough.

The paper is organized as follows. In Section \ref{s:mr} we introduce
hypotheses and notations and state main theorems. In Section
\ref{s:far} we construct an auxiliary function which represents the
principle part of the traveling wave solution for large $z$. Section
\ref{s:vs} describes the variational setting for the modified
functional.  In Section \ref{s:exist} we present results on existence
of traveling fronts and describe their properties.  Section
\ref{sec:uniq-supercr-fronts} is devoted to the proof of uniqueness of
supercritical fronts traveling with sufficiently large speed.  We also
discuss a simple one-dimensional example where the range of velocities
for which uniqueness holds can be explicitly computed.  Finally, in
Section \ref{sec:glob-asympt-stab} we prove a global stability result
for the supercritical fronts considered in Section
\ref{sec:uniq-supercr-fronts}.

\section{Preliminaries and main results}\label{s:mr}

In this section, we introduce the assumptions of our analysis and
state our main results. Throughout this paper we assume $\Omega$ to be
a bounded domain (connected open set, not necessarily simply
connected) with a boundary of class $C^2$. We start by listing the
assumptions on the nonlinearity $f$ which we will need in our
analysis:
\begin{enumerate}
\item[(H1)] The function $f: [0,1]\times\overline \Om \to \mathbb R$ 
  satisfies
  \begin{eqnarray} \label{c11} f(0,y) = 0, \qquad f(1, y) \leq
    0. \qquad \forall y \in \overline \Omega.
\end{eqnarray}
\item[(H2)] For some $\gamma \in (0,1)$
  \begin{eqnarray} \label{c12} f\in C^{0,\gamma}([0,1]\times\overline
    \Om), \qquad f_u\in C^{0,\gamma}([0,1]\times\overline \Om), 
\end{eqnarray}
where $f_u := \partial f / \partial u$.
\end{enumerate}

\noindent Note that hypotheses (H1) and (H2) are, in some sense, the
minimal assumptions needed to guarantee global existence and basic
regularity of solutions of \eqref{pde} satisfying
\begin{eqnarray}
  \label{u01}
  0 \leq u \leq 1.
\end{eqnarray}
From now on, when we speak of the solutions of either \eqref{pde} or
\eqref{trav}, we always assume that they satisfy \eqref{u01}.

Under the above hypotheses the following functional
\begin{eqnarray}
  \label{eq:Phic0}
  \Phi_c[u] := \int_\Sigma e^{cz} \left( \frac{1}{2} |\nabla
    u|^2 + V(u, y) \right) \, dx, \qquad V(u, y) := -\int_0^u f(s, y)
  \, ds, 
\end{eqnarray}
is well-defined for all functions obeying (\ref{bc}) and (\ref{u01})
that lie in the exponentially weighted Sobolev space $H^1_c(\Sigma)$,
i.e., functions which satisfy
\begin{eqnarray}
  \int_\Sigma e^{cz} (\left |\nabla u|^2 + u^2 \right) dx < \infty,
\end{eqnarray}
for a fixed $c > 0$.

Our final assumption on $f$ concerns the stability of the trivial
equilibrium $u = 0$. In this paper, we consider the situation, in
which this equilibrium is linearly {\em unstable}. This fact is
expressed by the following hypothesis:
\begin{enumerate}
\item[(U)] There holds
\begin{eqnarray}
  \label{eq:nu0}
  \nu_0 := \min_{\substack{\psi \in H^1(\Omega) \\ \psi|_{\partial
        \Omega_\pm} = 0}} \frac{\int_\Omega  ( |\nabla_y \psi|^2 -
    f_u(0, y) \psi^2 ) \, dy}{\int_\Omega  \psi^2 \, dy} < 0.
\end{eqnarray}
\end{enumerate}
Indeed, the constant $\nu_0$ defined in (\ref{eq:nu0}) is the
principal eigenvalue of the problem
\begin{eqnarray}
  \label{eq:psik}
  \Delta_y \psi_k + f_u(0, y) \psi_k + \nu_k \psi_k = 0,
\end{eqnarray}
with boundary conditions as in (\ref{bc}).  In particular, for the
principal eigenvalue we can assume $\psi_0$ to satisfy
\begin{eqnarray}
  \label{eq:psi0}
  \Delta_y \psi_0 + f_u(0, y) \psi_0 + \nu_0 \psi_0 = 0, \qquad \psi_0 >
  0, \qquad ||\psi_0||_{L^\infty(\Omega)} = 1.
\end{eqnarray}

Note that under hypothesis (U) we have existence of a particular
solution of (\ref{trav}) and (\ref{bc}) for some $c^* \geq c_0$, where
\begin{eqnarray}
  \label{eq:c0}
  c_0 := 2 \sqrt{-\nu_0} > 0,
\end{eqnarray}
which plays an important r\^ole for propagation with front-like initial
data \cite{mn:cms08}. As was shown in \cite{mn:cms08}, there are
precisely two alternatives: (i) either there exists a non-trivial
minimizer of $\Phi_{c^\dag}$ in $H^1_{c^\dag}(\Sigma)$ subject to
\eqref{bc} for some $c^\dag > c_0$, in which case $c^* = c^\dag$; or
(ii) there are no non-trivial minimizers of $\Phi_c$ in $H^1_c(\Sigma)$
subject to \eqref{bc} for all $c > c_0$, in which case $c^* = c_0$ and
there is a minimal wave with speed $c_0$. For our purposes here, it is
convenient to introduce the following equivalent characterization of
$c^*$:
\begin{eqnarray}
  \label{eq:cstar}
  c^* := \inf \left\{c > 0 : \inf_{0 \leq u \leq 1} \Phi_c[u] = 0 \right\}, 
\end{eqnarray}
where the infimum of $\Phi_c$ is taken over all functions in
$H^1_c(\Sigma)$ satisfying (\ref{bc}). The equivalence of the two
definitions of $c^*$ is shown by Proposition \ref{p:cstar}.

For $c > c^*$ any positive solution of (\ref{trav}) and (\ref{bc}) is
expected to satisfy
\begin{eqnarray}
  \label{eq:lammin}
  \bar u(y, z) \simeq a \psi_0(y) e^{-\lambda_-(c, \nu_0) z}, 
  \qquad   \lambda_\pm(c, \nu) = {c \pm \sqrt{c^2 + 4 \nu} \over 2},
\end{eqnarray}
as $z \to +\infty$, with some $a > 0$ and with $\psi_0 > 0$ being the
minimizer of the Rayleigh quotient in (\ref{eq:nu0}) normalized as in
(\ref{eq:psi0}). At least formally, the asymptotic behavior of
positive solutions of \eqref{trav} and \eqref{bc} can be obtained by
linearizing \eqref{trav} around $u = 0$, leading to the formula in
\eqref{eq:lammin} with the exponential decay rate determined by either
$\lambda_+(c, \nu_0)$ or $\lambda_-(c, \nu_0)$. However, the case
corresponding to $\lambda_+(c, \nu_0)$ is impossible, since then the
solution would belong to $H^1_c(\Sigma)$, i.e., it would be a
variational traveling wave (for the definition and further discussion,
see \cite{m:dcdsb04,lmn:arma08,mn:cms08}). But by \cite[Proposition
3.5]{lmn:arma08} this contradicts the assumption $c > c^*$. In other
words, for $c > c^*$ the traveling wave solutions acquire {\em fat
  tails}.

It is also well-known that solutions of \eqref{trav} become
$z$-independent as $z \to -\infty$. Thus, as $z \to -\infty$ the
function $\bar u(\cdot, z)$ is expected to approach a critical point
of the energy functional
\begin{eqnarray}
  \label{eq:E}
  E[v] = \int_\Omega \biggl( \frac12 |\nabla_y v|^2 +
  V(v, y) \biggr) 
  \, dy,
\end{eqnarray}
which is defined for all $v \in H^1(\Omega)$ with values in $[0, 1]$
satisfying Dirichlet boundary conditions on $\partial
\Omega_\pm$.\footnote{In the variational context throughout this
  paper, the Dirichlet boundary conditions are understood in the usual
  sense of considering the completion of the family of smooth
  functions whose support vanishes at the Dirichlet portion of the
  boundary with respect to the corresponding Sobolev norm.} By
hypothesis (H2) the critical points of $E$ satisfy the equation
\begin{eqnarray}
  \label{eq:v}
  \Delta_y v  + f(v, y) = 0, \qquad
  v|_{\partial \Omega_\pm} = 0, ~~ 
  \nu \cdot \nabla v |_{\partial \Omega_0} = 0.
\end{eqnarray}
Similarly to the case of $z \to +\infty$ discussed in the previous
paragraph, the convergence of the traveling wave solution to $v =
\lim_{z \to -\infty} \bar u(\cdot, z)$ satisfying \eqref{eq:v} can be,
at least formally, obtained from the analysis of \eqref{trav}
linearized around $v$. Introduce
\begin{eqnarray}
  \label{eq:tnu0}
  \tilde\nu_0 = \min_{\substack{\tilde\psi \in H^1(\Omega) \\
      \tilde\psi|_{\partial \Omega_\pm} = 0}} \frac{\int_\Omega  (
    |\nabla_y \tilde\psi|^2 - f_u(v, y) \tilde\psi^2 ) \,
    dy}{\int_\Omega  \tilde\psi^2 \, dy},
\end{eqnarray}
whose principal eigenvalue satisfies
\begin{eqnarray}
  \label{eq:tpsi0}
  \Delta_y \tilde\psi_0 + f_u(0, y) \tilde\psi_0 + \tilde\nu_0
  \tilde\psi_0 = 0, \qquad \tilde\psi_0 > 
  0, \qquad ||\tilde\psi_0||_{L^\infty(\Omega)} = 1.
\end{eqnarray}
Then in the generic case $\tilde \nu_0 > 0$, i.e., when $v$ is a
strict local minimizer of $E$, one expects the solution to approach
$v$ from below when $z \to -\infty$ as
\begin{eqnarray}
  \label{eq:vvv}
  \bar u(y, z) \simeq v(y) + \tilde a e^{|\lambda_-(c, \tilde \nu_0)| z}
  \tilde \psi_0(y) \qquad \text{for some} \qquad \tilde a < 0,
\end{eqnarray}
where we noted that the rate of convergence is determined by the only
negative value of $\lambda_\pm(c, \tilde\nu_0)$ \cite{mn:cms08}.

We now turn to our result on existence of solutions to \eqref{trav}
and \eqref{bc} that admit a particular variational
characterization. These are the traveling wave solutions which are
minimizers of the functional $\Phi_c^w$ defined as
\begin{eqnarray}
  \label{eq:Phic}
  \Phi_c^w[u] = \int_\Sigma e^{cz} \left( \frac{1}{2} |\nabla
    u|^2 + V(u, y) - \frac12 |\nabla w|^2 - V(w, y) \right) \, dx, 
\end{eqnarray}
for some suitably chosen $w: \Sigma \to \mathbb R$. Specifically, the
function $w \in C^2(\Sigma) \cap C^1(\overline\Sigma)$ should solve
\eqref{trav} and \eqref{bc} far ahead and far behind the front. As we
show in Sec. \ref{s:vs}, the functional $\Phi_c^w[u]$ is then well
defined for all $u$ satisfying \eqref{u01} and $u - w \in
H^1_c(\Sigma)$.

Our existence and variational characterization result is given by the
following theorem.

\begin{Theorem}
  \label{exist}
  Assume hypotheses (H1), (H2) and (U) are true. Then, for every $c >
  c^*$ and $a \in \mathbb R$ there exists $\bar u \in C^2(\Sigma) \cap
  C^1(\overline \Sigma) \cap W^{1,\infty}(\Sigma)$ such that $0 < \bar
  u < 1$ and:
  \begin{enumerate}
  \item[(i)] $\bar u$ solves (\ref{trav}) and (\ref{bc}), and there
    exists $w \in C^2(\Sigma) \cap C^1(\overline\Sigma)$ with values
    in $[0, 1]$ solving (\ref{trav}) and (\ref{bc}) in $\Omega \times
    (-\infty, -z_0) \cup \Omega \times (z_0, +\infty)$ for some $z_0 >
    0$, such that the function $\bar u$ minimizes $\Phi_c^w[u]$ over
    all $u$ satisfying \eqref{u01}, such that $u - w \in
    H^1_c(\Sigma)$.

  \item[(ii)] $\bar u(y, z) = a \psi_0(y) e^{-\lambda_-(c, \nu_0) z} +
    O(e^{-\lambda z})$, for some $\lambda > \lambda_-(c, \nu_0)$ and
    $\psi_0$ solving (\ref{eq:psi0}), uniformly in
    $C^1(\overline\Omega \times [R, +\infty))$, as $R \to +\infty$.

  \item[(iii)] $\bar u$ is a strictly monotonically decreasing
    function of $z$, and we have $\bar u(\cdot, z) \to v$ uniformly in
    $C^1(\overline\Omega)$, where $0 < v \leq 1$ is a critical point
    of $E$, and $E[v] < 0$.

  \item[(iv)] $\tilde \nu_0 \ge 0$, moreover, if $\tilde \nu_0 > 0$,
    then $\bar u(y,z) = v(y) + \tilde a \tilde \psi_0(y)
    e^{-\lambda_-(c, \tilde \nu_0)z} + O(e^{-\lambda z})$, with
    some $\tilde a < 0$ and $\lambda < \lambda_-(c, \tilde
    \nu_0)$, uniformly in $C^1(\overline\Omega \times (-\infty, R])$,
    as $R \to -\infty$.

  \item[(v)] If $\tilde w \in C^2(\Sigma) \cap C^1(\overline \Sigma)$
    with values in $[0, 1]$ solves (\ref{trav}) and (\ref{bc}) in
    $\Omega \times (-\infty, -\tilde z_0) \cup \Omega \times (\tilde
    z_0, +\infty)$ for some $\tilde z_0 > 0$ and $\tilde w - w \in
    H^1_c(\Sigma)$, then $\bar u$ also minimizes $\Phi_c^{\tilde w}$.
  \end{enumerate}
\end{Theorem}

Let us point out that existence results for supercritical fronts in
cylinders were first obtained in the classical work of Berestycki and
Nirenberg under a few extra assumptions \cite{berestycki92} (see also
\cite{vega93,vega93cpde}). Specifically, they assumed (in our notation
and setting) that if $v_1$ is the smallest positive critical point of
$E$ (in the sense that if $0 < v \leq 1$ is a critical point of $E$,
then $v \geq v_1$) and if $v_1$ is stable, i.e., if $\tilde \nu_0 > 0$
with $v = v_1$ in \eqref{eq:tnu0}, then there exists a unique (up to
translations) traveling wave solution approaching zero from above as
$z \to +\infty$ and $v_1$ from below as $z \to -\infty$. Note that
existence of $v_1$ and the property $\tilde\nu_0 \geq 0$ for $v = v_1$
follows immediately from our hypotheses.

It is easy to see that the existence and monotonicity result for
fronts connecting zero and $v_1$ may also be obtained from a variant
of our Theorem \ref{exist}, in which a minimizer of $\Phi_c^w$ is
constructed in the class of functions bounded above by $v_1$. More
precisely, let us define the value of $c_1^*$ by \eqref{eq:cstar}, in
which the infimum is taken over all $u \in H^1_c(\Sigma)$ satisfying
\eqref{bc} and such that $u \leq v_1$:
\begin{eqnarray}
  \label{eq:26}
  c_1^* = \inf \left \{ c > 0 : \inf_{0 \leq u \leq v_1} \Phi_c[u] = 0
  \right\}, 
\end{eqnarray}
where, again, the infimum is taken over all $u \in H^1_c(\Sigma)$
satisfying \eqref{bc}.  Then under hypotheses (H1), (H2), and (U), for
every $c > c_1^* \geq c_0$ and every $a \in \mathbb R$ there exists
$\bar u$ which satisfies the conclusions of Theorem \ref{exist} with
$v = v_1$, where $v_1$ is the smallest positive solution of
\eqref{eq:v}. In particular, it can be seen from Theorem \ref{exist}
that $\bar u$ is the minimizer of $\Phi_c^w$ for all $w$ with values
in $[0, 1]$, such that $w - \bar u \in H^1_c(\Sigma)$.

It is important to point out that using their sliding domains method,
Berestycki and Nirenberg also proved uniqueness of solutions of
\eqref{trav} and \eqref{bc} connecting zero and a non-degenerate $v_1$
for all $c > c_1^*$ (see also \cite{vega93,vega93jmma,vega93cpde}). On
the other hand, a priori one cannot expect uniqueness in the more
general class of solutions satisfying \eqref{u01}, even in one
dimension (see the counter-example in
Sec. \ref{sec:uniq-supercr-fronts}).  A natural question, then, is
whether for a given $c > c^*$ the minimizer given by Theorem
\ref{exist} is the solution obtained by Berestycki and Nirenberg. We
can give a negative answer to this question in the case when $c^* >
c_1^*$ for a certain range of $c > c^*$.

Indeed, suppose $c^* > c_1^*$. Then, since $c_1^* \geq c_0$, there
exists a non-trivial minimizer $\bar u^*$ of $\Phi_{c^*}$ among all $0
\leq u \leq 1$. Furthermore, since $\bar u^*$ is not a minimizer of
$\Phi_{c^*}[u]$ over all $u \leq v_1$, we have $\lim_{z \to -\infty}
\bar u^*(\cdot, z) > v_1$. In particular, this means that for
\begin{eqnarray}
  \label{eq:cdv}
  c^\dag_{v_1} :=  \inf \left \{ c > 0 : \inf_{v_1 \leq u \leq 1}
    \Phi_c^{v_1}[u] = 0 \right\}
\end{eqnarray}
we have $c^\dag_{v_1} \geq c^*$. If not, for $u_1 = \min (\bar u^*,
v_1)$ and $u_2 = \max(\bar u^*, v_1)$ we would obtain $\Phi_{c^*}[\bar
u^*] = \Phi_{c^*}[u_1] + \Phi_c^{v_1}[u_2] \geq \Phi_{c^*}[u_1]$,
implying that $u_1 \leq v_1$ is a non-trivial minimizer. However, this
contradicts the fact that there are no non-trivial minimizers of
$\Phi_{c^*}$ that lie below $v_1$.

In fact $c^\dag_{v_1} > c^*$, since the infimum in \eqref{eq:cdv} is
attained, with the value of $c^\dag_{v_1}$ being the speed of the
non-trivial minimizer of $\Phi_{c^\dag_{v_1}}^{v_1}$ over all $u$,
such that $u - v_1 \in H^1_c(\Sigma)$ and $v_1 \leq u \leq 1$. Indeed,
there exists a sequence $c_n \to c^\dag_{v_1}$ from below and a
sequence of functions $u_n \geq v_1$, such that $\Phi_{c_n}^{v_1}[u_n]
< 0$. Therefore, in view of the fact that $\tilde \nu_0 \geq 0$ for $v
= v_1$, by \cite[Theorem 3.3]{mn:cms08} applied to $\Phi_c^{v_1}$ (see
also the proof of \cite[Lemma 3.5]{mn:sima11}), there exists a
non-trivial minimizer of $\Phi_c^{v_1}$ for some $c \geq
c^\dag_{v_1}$, and by the definition of $c^\dag_{v_1}$ and the
arguments in the proof of Proposition \ref{p:cstar} we have $c =
c^\dag_{v_1}$.

Since by the argument in the proof of \cite[Lemma 3.5]{mn:sima11} for
every $0 \leq u \leq v_1$, such that $u - w \in H^1_c(\Sigma)$, with
$w$ as in Theorem \ref{exist}, and $u(\cdot, z) \to v_1$ in
$C^1(\overline\Omega)$ as $z \to -\infty$, the value of $\Phi_c^w[u]$
can be lowered by extending $u$ above $v_1$ for every $c <
c^\dag_{v_1}$, we have the following consequence of Theorem
\ref{exist}:

\begin{Corollary}
  \label{c:2}
  Let $c^* > c_1^*$. Then $c^\dag_{v_1} > c^*$, and for every $c^* < c
  < c^\dag_{v_1}$ the solution $\bar u$ in Theorem \ref{exist} has $v
  > v_1$.
\end{Corollary}

\noindent In other words, if $c^* > c_1^*$, for every $c$ as in
Corollary \ref{c:2} and every $a \in \mathbb R$ there exist at least
two solutions of \eqref{trav} and \eqref{bc} satisfying
\eqref{eq:lammin}. 

\begin{Remark}
  \label{r:2}
  Note that the argument leading to Corollary \ref{c:2} also works
  under the assumption that $c^\dag_{v_1} > c^*$ when $c^* = c_1^* =
  c_0$.
\end{Remark}

Let us now go back to the question under which conditions the
solutions of Berestycki and Nirenberg, which connect zero with $v_1$
for every $c > c_1^*$ are, in fact, the unique (up to translations)
solutions of \eqref{trav} and \eqref{bc}. In one dimension it is
possible to use phase plane arguments to show that for large enough
values of $c$ the only solution of \eqref{trav} connects zero to $v_1$
(e.g., under non-degeneracy of $v_1$). What about the
higher-dimensional case? In the following we show that the solution
connecting zero to $v_1$ is the only traveling wave solution for large
enough speeds, up to translations.

We start by introducing
\begin{eqnarray}
  \label{eq:hnu0}
  \hat\nu_0 := \min_{\substack{\psi \in H^1(\Omega) \\
      \psi|_{\partial \Omega_\pm} = 0}} \frac{\int_\Omega  (
    |\nabla_y \psi|^2 - q(y) \psi^2 ) \,
    dy}{\int_\Omega  \psi^2 \, dy}, \qquad q(y) := \max_{s \in [0, 1]}
  f_u(s, y). 
\end{eqnarray} 
In analogy with \eqref{eq:c0}, we can then define the quantity
\begin{eqnarray}
  \label{eq:c1}
  c^\sharp := 2 \sqrt{-\hat\nu_0},
\end{eqnarray}
which is easily seen to satisfy $c_0 \leq c^\sharp < \infty$. It turns
out that $c^\sharp$ gives a lower bound for the values of $c$, for
which the solution of \eqref{trav} and \eqref{bc} is unique, up to
translations. Our global uniqueness result, which relies on the
variational characterization of the obtained solutions is contained in
the following theorem.

\begin{Theorem}
  \label{unique}
  Let $c > c^\sharp$, where $c^\sharp \geq c_0$ is defined in
  \eqref{eq:c1}, and let the assumptions of Theorem \ref{exist} be
  satisfied. Then $c^\sharp \geq c^*$, and for every $c > c^\sharp$
  there exists a unique solution of \eqref{trav} and \eqref{bc} (up to
  translations).
\end{Theorem}

\noindent Furthermore, the unique solution of Theorem \ref{unique} is
globally stable with respect to perturbations with sufficiently fast
exponential decay.

\begin{Theorem}
  \label{converge}
  Let $\bar u$ be the traveling wave solution with speed $c >
  c^\sharp$ from Theorem \ref{unique}, and let $0 \leq u_0 \leq 1$
  satisfy $u_0 -\bar u \in L^2_{c'}(\Sigma)$, for some $c' > 2
  \lambda_-(c, \hat\nu_0)$. Then, if $u$ is the solution of
  (\ref{pde}) and (\ref{bc}) with $u(x, 0) = u_0(x)$, we have
  \begin{eqnarray}
    ||u(y, z + c t, t) - \bar u(y, z)||_{L^2_{c'}(\Sigma)} \leq C
    e^{-\sigma t}, 
  \end{eqnarray}
  for some $\sigma > 0$, $C > 0$.
\end{Theorem}

\noindent We note that the result of Theorem \ref{converge} also holds
in other norms, e.g. with respect to uniform convergence on compacts
(see \cite{mn:sima11} for details).

Note that as a direct consequence of Theorems
\ref{exist}--\ref{converge} we have a complete characterization of the
supercritical fronts in problems involving nonlinearities of the kind
that was originally considered in \cite{kolmogorov37,fisher37}:
\begin{Corollary}
  \label{c:kpp}
  Under hypotheses (H1), (H2) and (U), assume that, in addition,
  $f_u(u, \cdot) \leq f_u(0, \cdot)$ for all $u \in [0, 1]$. Then for
  each $c > c_0$ there exists a unique (up to translations) solution
  of \eqref{trav} and \eqref{bc}, which is globally stable in the
  sense of Theorem \ref{converge} with respect to perturbations in
  $L^2_{c'}(\Sigma)$, with any $c' > 2 \lambda_-(c, \nu_0)$.
\end{Corollary}
\noindent This corollary is obtained by noting that under its
assumption on the nonlinearity we simply have $\nu_0 = \hat
\nu_0$. Therefore, $c^\sharp = c_0$ and, in particular, $c^* = c_0$ as
well. In fact, in view of the conclusion of Theorem \ref{exist}(ii)
the following convergence result holds under a single assumption on
the rate of exponential decay of the initial data.

\begin{Corollary}
  \label{c:kppexp}
  Under the assumptions of Corollary \ref{c:kpp}, let $c > c_0$ and
  assume that $u_0(y, z) = a \psi_0(y) e^{-\lambda_-(c, \nu_0) z} +
  O(e^{-\lambda z})$, for some $a \in \mathbb R$ and $\lambda >
  \lambda_-(c, \nu_0)$. Then the conclusion of Theorem \ref{converge}
  holds.
\end{Corollary}

In the rest of the paper (H1), (H2), (U) are always assumed to be
satisfied.

\section{The front far ahead}
\label{s:far}

In this section we start the investigation of supercritical traveling
wave solutions.  As was mentioned in the introduction, the first
important step in the study of traveling fronts is the construction of
an auxiliary function which solves (\ref{trav}) and \eqref{bc} for $z$
large positive and has the asymptotic behavior given by
(\ref{eq:lammin}). Specifically, we will construct a function $w$ with
the following properties:
\begin{enumerate}
\item[(W1)] $w \in C^2(\Sigma) \cap C^1(\overline \Sigma)$ and solves
\begin{eqnarray}
  \label{eq:w}
  \Delta w + c w_z +f(w, y) = 0, 
\end{eqnarray}
in $\Sigma$ for $z > 3$, with the same boundary conditions as in
(\ref{bc}). 

\item[(W2)] There exists $\delta > 0$ and $a > 0$ sufficiently small,
  such that
\begin{eqnarray}
  \label{eq:w0}
  w(y, z) = a \psi_0(y) e^{-\lambda_-(c, \nu_0) z}(1 + O(e^{-\delta
    z})), \qquad z \to +\infty,
\end{eqnarray}
where $\lambda_-(c, \nu_0)$ and $\psi_0$ are defined in
\eqref{eq:lammin} and \eqref{eq:psi0}.

\item[(W3)] We have
  \begin{eqnarray}
    \label{eq:w2a}
    0 \leq w(y, z) \leq 2 a \psi_0(y) \leq 1 \quad \forall (y,
    z) \in \Sigma, 
  \end{eqnarray}
  and $w = 0$ for all $z < 2$.
\end{enumerate} 

In order to establish existence of a function $w$ satisfying (W1) --
(W3), we only need to prove existence of a solution of (\ref{eq:w})
satisfying (W2) and (W3) for $z > 2$. Then we obtain the desired
function by multiplying this solution with a cutoff function $\eta \in
C^\infty(\mathbb R)$, such that $\eta(z) = 1$ for $z > 3$ and $\eta(z)
= 0$ for $z < 2$.

By translational symmetry in the $z$-direction, any translate of
solution of (\ref{eq:w}) is once again a solution for sufficiently
large $z$.  Translations of $w$ are equivalent to adjusting the value
of the constant $a$ in (\ref{eq:w0}). Therefore, fixing the constant
$a$ in (\ref{eq:w0}) is equivalent to fixing translations of the
solution of (\ref{eq:w}).  Thus, the choice of a particular value of
$a$ is inconsequential. Indeed, considering solutions with $a$ fixed
and $z$ large enough positive is equivalent to considering the problem
for small enough $a>0$ and $z>0$.  Conceptually, this choice of $w$
resembles selecting one element of the stable manifold corresponding
to the fixed point $u=0$ in the one-dimensional setting (in the sense
of ordinary differential equations).

The following proposition establishes existence of a function $w$ with
all the desired properties.


\begin{Proposition}
  \label{p:wexists}
  For every $c > c_0$ and every $a > 0$ sufficiently small, there
  exists a function $w$ with properties (W1)--(W3).
\end{Proposition}

The proof of this proposition is based on two lemmas.  We first
explicitly construct sub and supersolutions (see definitions (i) and
(ii) below) for the problem (\ref{eq:w}) on the right half cylinder
(Lemma \ref{l:subsup}).  Next we show that if one restrict the problem
(\ref{eq:w}) on bounded sections of the cylinder, then this problem
has a solution squeezed between sub and supersolutions (Lemma
\ref{l:wexR}). We then complete the proof of Proposition
\ref{p:wexists} by passing to the limit.

In order to proceed, we introduce the notions of sub and supersolutions:
\begin{enumerate}
\item[(i)] We call $\bw$ a supersolution of the problem (\ref{eq:w}),
  if $\bw \in C^2(\Omega \times (0, +\infty))$ satisfies the boundary
  conditions in (\ref{bc}) and
\begin{eqnarray}\label{bw}
\Delta \bw+c\bw_z+f(\bw,y)\le 0, \qquad z>0,
\end{eqnarray}
with
\begin{eqnarray} \label{bw0}
  \bw(0,y) = 2 a \psi_0(y), \qquad y\in \Omega.
\end{eqnarray}

\item[(ii)] We call $\uw$ a subsolution for the problem (\ref{eq:w}),
  if $\uw \in C^2(\Omega \times (0, +\infty))$ satisfies the boundary
  conditions in (\ref{bc}) and
\begin{eqnarray}\label{uw}
\Delta \uw+c\uw_z+f(\uw,y)\ge 0, \qquad z>0,
\end{eqnarray}
with 
\begin{eqnarray}\label{uw0}
\uw(0,y) = 0, \qquad y\in \Omega.
\end{eqnarray}
\end{enumerate}

\begin{Lemma}
  \label{l:subsup}
  There exist $\underline{w}$ and $\overline{w}$, the sub- and
  super-solutions of (\ref{eq:w}), which satisfy (W2) and (W3), for
  every $a > 0$ small enough.
\end{Lemma}

\begin{proof} The proof follows by explicit construction of $\uw$ and
  $\bw$ verifying (\ref{bc}), (\ref{bw}) -- (\ref{uw0}), (W2), and
  (W3). By assumption (H2), for $||w||_{L^\infty(\Sigma)} \leq 2 a \ll
  1$ there exist constants $\uk, \bk > 0$, such that
  \begin{eqnarray}
    f_u(0,y)w-\uk |w|^{1+\gamma}\le f(w,y)\le f_u(0,y)w+\bk
    |w|^{1+\gamma}.  
\end{eqnarray}
Therefore, the function $\bw$ satisfying (\ref{bc}), (\ref{bw0}) and
\begin{eqnarray}
  \label{Nw}
  \mathcal N(\bar w) := \Delta \bw +c \bw_z+f_u(0,y)\bw+\bk
  \bw^{1+\gamma} \le 0,
\end{eqnarray}
in $\Omega \times (0, +\infty)$ automatically satisfies (\ref{bw}),
and thus is a supersolution.

Let us show that the function 
\begin{eqnarray}
  \label{bww}
  \bw(z,y)=a\psi_0(y)e^{-\lambda_- z}(1 + e^{-\delta z})
\end{eqnarray}
is a supersolution for some $0 < \delta \ll 1$ (throughout the proof
we use the shorthand $\lambda_-$ for $\lambda_-(c, \nu_0)$). Indeed,
substituting this function into (\ref{Nw}), we obtain
\begin{eqnarray}
  \mathcal N(\bw) & = & 
  a ( \Delta_y \psi_0 + (\lambda_-^2 - c \lambda_-
  + f_u(0, y)) \psi_0 ) 
  e^{-\lambda_- z}(1 + e^{-\delta z}) 
  \nonumber \\ && - a \delta ( c - 2 \lambda_- - \delta) \psi_0
  e^{-(\lambda_- + \delta) z}   \label{Nww} \\ && + \bk a^{1+\gamma}
  \psi_0^{1+\gamma} 
  e^{-(1 + \gamma) \lambda_- z} (1 + e^{-\delta
    z})^{1+\gamma}. \nonumber 
\end{eqnarray}
The first line in (\ref{Nww}) is identically zero in view of the
definition of $\lambda_-$ and $\psi_0$, see (\ref{eq:psi0}) and
(\ref{eq:lammin}). Therefore, 
\begin{eqnarray}
  \label{eq:Nwww}
  \mathcal N(\bw) & = & - a \delta ( c - 2 \lambda_- - \delta) \psi_0 
  e^{-(\lambda_- + \delta) z} \nonumber \\ && \times \left\{ 1 -
    \frac{\bar k a^\gamma \psi_0^\gamma  (1 + e^{-\delta
        z})^{1 + \gamma}}{\delta (c - 2 \lambda_- - \delta)} e^{-(\gamma
      \lambda_- - \delta) z} \right\}. 
\end{eqnarray}
Since also $c > 2 \lambda_-$ by (\ref{eq:lammin}), the factor
multiplying the expression in the curly brackets in (\ref{eq:Nwww}) is
negative for $\delta$ sufficiently small. Let us show that the
expression in the curly brackets can be made positive by a suitable
choice of $a$ and $\delta$. Indeed, since the expression in the curly
bracket is an increasing function of $z$ for $\delta$ small enough, we
only need to verify its positivity at $z = 0$. Then, the result easily
follows by choosing $\delta$ sufficiently small first, and then
choosing $a$ such that $a^\gamma \ll \delta$.

We now turn to sub-solutions. By the same argument as above, it is
enough to show that
\begin{eqnarray}
  \label{Mw}
  \mathcal M(\uw) := \Delta \uw +c \uw_z+f_u(0,y)\uw-\uk
  \uw^{1+\gamma} \ge 0
\end{eqnarray}
in $\Omega \times (0, +\infty)$ for some $\uw \geq 0$. Similarly to
(\ref{bww}), we choose
\begin{eqnarray}
  \label{uww}
  \uw = a \psi_0 e^{-\lambda_- z} (1 - e^{-\delta z}).
\end{eqnarray}
Substituting this into (\ref{Mw}) and using the same arguments as for
the supersolution, we obtain
\begin{eqnarray}
  \label{eq:Mwww}
  \mathcal M(\uw) & = & a \delta ( c - 2 \lambda_- - \delta) \psi_0 
  e^{-(\lambda_- + \delta) z} \nonumber \\ && \times \left\{ 1 -
    \frac{\uk a^\gamma \psi_0^\gamma  (1 - e^{-\delta
      z})^{1 + \gamma}}{\delta (c - 2 \lambda_- - \delta)} e^{-(\gamma
    \lambda_- - \delta) z} \right\}, 
\end{eqnarray}
which is positive for a similar choice of $\delta$ and $a$.

Finally, we note that, as a consequence of our construction, the
functions $\bw$ and $\uw$ satisfy (\ref{eq:w0}) and (\ref{eq:w2a}).
\end{proof}

We next establish existence of a classical solution of (\ref{eq:w}) on
finite sections of the cylinder away from the boundaries. 

\begin{Lemma}
  \label{l:wexR}
  Let $R > 1$. Then there exists $w_R \in C^2(\Omega \times (1, R)) \cap
  C^1(\overline \Omega \times [1, R])$ solving (\ref{eq:w}) in
  $\overline \Omega \times [1, R]$, and satisfying $\uw \leq w_R \leq
  \bw$, where $\uw$ and $\bw$ are as in Lemma \ref{l:subsup}.
\end{Lemma}

\begin{proof}
  First, the proof of existence of a weak solution of (\ref{eq:w}) in
  $\Sigma_R = \Omega \times (0, R+1)$, sandwiched between the sub- and
  the super-solutions of Lemma \ref{l:subsup} follows by standard
  monotone iteration argument (see, e.g., \cite[Sec. 9.3]{evans}.)
  Furthermore, by standard elliptic regularity, this solution is
  regular in the interior of $\Sigma_R$ \cite{gilbarg}. Now, we apply
  $L^p$ regularity theory to the solution in $\widetilde \Sigma_R =
  \Omega \times (1, R) \subset \subset \Sigma_R$ \cite[Theorem
  9.13]{gilbarg} to obtain regularity up to the boundary on the
  Dirichlet portions $\partial \Sigma_\pm$ of the cylinder:
  \begin{eqnarray}
    \label{eq:wwp}
    ||w_R||_{W^{2,p}(\widetilde \Sigma_R)} \leq C (||w_R||_{L^p(\Sigma_R)}
    + ||f||_{L^p(\Sigma_R)}) \leq C' ||w_R||_{L^p(\Sigma_R)},
  \end{eqnarray}
  for some $C, C' > 0$ and any $p < \infty$, where we took into
  account hypothesis (H2) on $f$. Hence, by Sobolev imbedding $w_R$ is
  in $C^{1,\alpha}(\widetilde\Sigma_R)$, for any $\alpha \in (0,
  1)$. By construction, $\uw \leq w_R \leq \bw$ in the closure of
  $\widetilde\Sigma_R$.
\end{proof}

\begin{proof}[Proof of Proposition  3.1]
  We now finish the proof of Proposition \ref{p:wexists} by passing to
  the limit $R \to \infty$ in Lemma \ref{l:wexR}. Observe that by the
  bounds of Lemma \ref{l:wexR} and the exponential decay of $\uw$ and
  $\bw$ of Lemma \ref{l:subsup}, we have $||w_R||_{L^p(\Omega \times
    (1, R))} \leq C$, where $C > 0$ is independent of $R$. Therefore,
  by the estimate of (\ref{eq:wwp}) we have $w_R \to w$ in
  $W^{2,p}_\mathrm{loc}(\Omega \times (1, \infty))$ on a sequence of
  $R_k \to \infty$. In view of arbitrariness of $p < \infty$, the
  convergence of $w_R$ to $w$ is also in $C^1(\overline \Omega \times
  [2, R-1])$, for any $R > 3$. Hence, the limit $w$ is a weak solution
  of (\ref{eq:w}) in $\overline \Omega \times [2, \infty)$, and by
  standard elliptic regularity is also a classical solution. The proof
  is completed by multiplying the obtained solution by a cutoff
  function.
\end{proof}

\begin{Remark} \label{r:walot}
  Our proof guarantees existence of infinitely many functions
  satisfying (W1) -- (W3), since existence of a solution for
  (\ref{eq:w}) also holds for each fixed boundary data at $z = 0$
  satisfying (\ref{eq:w2a}).
\end{Remark}

\section{Variational setting}
\label{s:vs}
As was already noted in the Introduction, certain solutions of
(\ref{trav}) may be viewed as critical points of the functional
$\Phi_c[u]$ introduced in (\ref{eq:Phic0}). This, however, poses some
restrictions on the type of traveling wave solutions (which are the
so-called variational traveling waves
\cite{m:dcdsb04,mn:cms08,lmn:arma08}) captured by this
approach. Indeed, in view of the presence of the exponential weight in
the functional the set of admissible functions must be characterized
by the natural exponential decay at $z = +\infty$. At the same time,
it is known (see, e.g. \cite{berestycki92,vega93}) that when $u = 0$
is an {\em unstable} equilibrium of (\ref{pde}), generically there are
traveling wave solutions which fail to have this natural exponential
decay. In this situation the use of a variational approach to the
construction of traveling wave solutions seems problematic due the
divergence of the respective integrals, even though formally
(\ref{trav}) is still the Euler-Lagrange equation for $\Phi_c$.  We
propose to eliminate this divergence by subtracting the divergent part
of the integrand in the definition of $\Phi_c$. Namely, we start by
introducing the modified functional in \eqref{eq:Phic}, where $w$ is a
function satisfying properties (W1) -- (W3) of Sec. \ref{s:far}. In
other words, we subtract from the integrand in the definition of
$\Phi_c$ its value obtained by substituting an exact solution of
(\ref{trav}) for $z > 3$ (property (W1) of $w$) and characterized by a
prescribed exponential decay at $z = +\infty$ (property (W2) of
$w$). Note that existence of $w$ with these properties was proved in
Sec. \ref{s:far} by non-variational techniques.

We next identify the suitable admissible class of functions $u$ on
which the new functional $\Phi_c^w[u]$ makes sense.  Let $u = w + h$,
for some unknown function $h : \Sigma \to \mathbb R$. With this
definition, we now have $\Phi_c^w[u] = \Psi_c^w[h]$, where
\begin{eqnarray}
  \label{eq:Phich}
  \Psi_c^w[h] = \int_\Sigma e^{cz} \biggl( \frac{1}{2} |\nabla
  h|^2 + V(w + h, y) - V(w, y) - V'(w, y) h \biggr) dx \nonumber \\ 
  - \int_\Sigma e^{cz} \biggl( \Delta w + c w_z - V'(w, y) \biggr) h
  \, dx,
\end{eqnarray}
where here and throughout the rest of this section we use the notation
$ V'(u,y) := \partial V(u,y) / \partial u$ and $V''(u,y) := \partial^2
V(u,y) / \partial u^2$.  Observe that as long as $u(x) \in [0,1]$ for
all $x \in \Sigma$, we can use Taylor formula to rewrite $\Psi_c^w[h]$
in the following form
\begin{eqnarray}
  \Psi_c^w[h] = \int_\Sigma e^{cz} \biggl( \frac{1}{2} |\nabla
  h|^2 + \frac12 V''(\tilde u, y) h^2 \biggr) dx \nonumber \\ 
  - \int_2^3 \int_\Omega e^{cz} \biggl( \Delta w + c w_z -
  V'(w, y) \biggr) h \, dy dz,
\end{eqnarray}
for some $\tilde u$ sandwiched between $w$ and $u$. Here we took into
account hypothesis (H2) and that, by property (W1) of $w$, the
integrand in the last term in (\ref{eq:Phich}) is identically zero for
$z > 3$ and $z < 2$.  Now, since $V', V'', w, w_z$ are uniformly
bounded and $\Delta w \in L^2(\Omega \times (2,3))$, the functional
$\Psi_c^w$ is well-defined within the admissible class $\mathcal
A_c^w$, defined as follows:
\begin{definition}
  A function $h: \Sigma \to \mathbb R$ is said to belong to the
  admissible class $\mathcal A_c^w$, if $h \in H^1_c(\Sigma)$ and $0
  \leq h(x) + w(x) \leq 1$ for all $x \in \Sigma$.
\end{definition}
Note that the assumption that $h \in \mathcal A_c^w$ implies that $u =
w + h$ has values in the unit interval is not too restrictive, in view
of hypothesis (H1), since $u = 0$ and $u = 1$ are sub- and
super-solutions of (\ref{trav}) and, hence, are the natural barriers
\cite{kinderlehrer}.

In order to apply the direct method of calculus of variation, we need
to establish weak sequential lower-semicontinuity and coercivity of
the functional $\Psi_c^w$ in the considered function
class. Lower-semicontinuity in the considered class of problems was
studied in \cite{lmn:arma08}, and we have
\begin{Proposition}
  \label{p:lsc}
  For all $c > c_0$, the functional $\Psi_c^w$ is sequentially
  lower-semicontinuous in the weak topology of $H^1_c(\Sigma)$ within
  $\mathcal A_c^w$.
\end{Proposition}
\begin{proof}
  The proof follows exactly as in \cite[Proposition 5.5]{lmn:arma08}.
\end{proof}

\noindent Thus, the main difficulty in applying a variational approach
to supercritical fronts has to do with coercivity of $\Psi_c^w$. Below
we show that coercivity indeed holds for all $c > c^*$.

\begin{Proposition}
  \label{p:coer}
  The functional $\Psi_c^w$ is coercive in $\mathcal A_c^w$ for every
  $c > c^*$.
\end{Proposition}

\begin{proof}
  We carry out the proof in 4 steps.

  \paragraph{Step 1}
  First of all, since $h, V', w, w_z$ are uniformly bounded and since
  $\Delta w \in L^2(\Omega \times (2,3))$, by a straightforward
  application of Cauchy-Schwarz inequality we have
  \begin{eqnarray}
    \left| \int_2^3\int_\Omega e^{cz} \biggl( \Delta w + c w_z - V'(w, y)
      \biggr) h \, dydz \right|  \leq K_1,
  \end{eqnarray}
  for some $K_1 > 0$. On the other hand, using the formula in the
  proof of \cite[Proposition 6.9]{lmn:arma08} and the fact that by the
  definition of $c^*$ we have $\Phi_{c'}[\rho] \geq 0$ for all $c' \in
  (c^*, c)$ and all $\rho \in H^1_{c'}(\Sigma)$, it holds that
  \begin{eqnarray}
    \Phi_c[h] = \int_\Sigma e^{cz} \biggl(
    \frac12 |\nabla h|^2 + V(h, y) \biggr) dx 
    \geq {c^2 - {c'}^2 \over 2 c^2} \int_\Sigma e^{cz} h_z^2 \,  dx.
  \end{eqnarray}
  Therefore, we have the following estimate:
\begin{eqnarray}
  \Psi_c^w[h]  + K_1 \geq  \int_\Sigma e^{cz} \biggl(
  \frac12 |\nabla h|^2 + V(h, y) \biggr) dx \nonumber \\ 
  + \int_\Sigma e^{cz} \biggl( V(w + h, y) - V'(w, y)
  h - V(w, y) - V(h, y) \biggr) dx \\ 
  \geq {c^2 - {c'}^2 \over 2 c^2} \int_\Sigma e^{cz} h_z^2 \,  dx
  \nonumber \\   
  + \int_\Sigma e^{cz} \biggl( V(w + h, y) - V'(w, y)
  h - V(w, y) - V(h, y) \biggr) dx. \nonumber
\end{eqnarray}

\paragraph{Step 2}

By Taylor formula, there exist $\tilde w$ and $\tilde h$ such that
$|\tilde w - w| < |h|$ and $|\tilde h| < |h|$ such that
\begin{eqnarray}
  V(w + h, y) - V'(w, y) h - V(w, y) - V(h, y) \nonumber \\ 
  = \frac12 \left( V''(\tilde w, y) - V''(\tilde h, y) \right)
  h^2. \label{eq4_6}
\end{eqnarray}
On the other hand, hypothesis (H2) implies that
\begin{eqnarray}
  \label{agam}
  |V''(\tilde w, y) - V''(\tilde h, y)| \leq C_1 |\tilde w - \tilde
  h|^\gamma,
\end{eqnarray}
for some $C_1 > 0$.  Therefore, for all $|h| \leq a$, where $a$ is
defined in property (W3) of $w$, we obtain
\begin{eqnarray}
  | V(w + h, y) - V'(w, y) h - V(w, y) - V(h, y) | \leq C_2
  a^\gamma h^2, \label{eq4_8}
\end{eqnarray}
for some $C_2 > 0$.

Now, by boundedness of $V', V''$, when $|h| > a$ we also have
\begin{eqnarray}
  \label{agaw}
  | V(w + h, y) - V'(w, y) h - V(w, y) - V(h, y) | \leq C_3 w,
\end{eqnarray}
for some $C_3 > 0$.

\paragraph{Step 3}

Using \eqref{eq4_6}--(\ref{eq4_8}) and the Poincar\'e inequality
\cite[Lemma 5.1]{lmn:arma08}, one can write
\begin{eqnarray}
  \int_{\{|h| \leq a \}}  e^{cz} \bigl| V(w + h, y) - V'(w, y)
  h - V(w, y) - V(h, y) \bigr| dx \nonumber \\
  \leq C_2 a^\gamma \int_\Sigma e^{cz} h^2 \, dx \leq {4 C_2
    a^\gamma \over c^2} \int_\Sigma e^{cz} h_z^2 \,  dx.
\end{eqnarray}

Similarly, using (\ref{agaw}) and property (W3) of $w$ we have
\begin{eqnarray}
  \int_{\{|h| > a \}}  e^{cz} \bigl| V(w + h, y) - V'(w, y)
  h - V(w, y) - V(h, y) \bigr| dx \nonumber \\
  \leq C_3 \int_{-\infty}^R\int_{\{|h(\cdot, z)| > a \}}  e^{cz} w \,
  dy dz +  C_3 \int_R^\infty \int_{\{|h(\cdot, z)| > a \}}  e^{cz} w \,
  dy dz \nonumber \\ \leq {2 a C_3 |\Omega| e^{c R} \over
    c} + C_3 \int_R^\infty \int_{\{|h(\cdot, z)| > a \}} e^{cz} w \,
  dy dz,
\end{eqnarray}
for any $R > 0$, where $|\Omega|$ denotes the $(n-1)$-dimensional
Lebesgue measure of $\Omega$. Then, since by Poincar\'e inequality
\cite[Lemma 2.1]{lmn:arma08} we have
\begin{eqnarray} \label{huliRz}
  {e^{-c z} \over c} \int_z^\infty \int_\Omega e^{c\zeta} h_{\zeta}^2
  \,  dy d\zeta \geq   
  \int_\Omega h^2(y, z) \, dy \nonumber \\ \geq 
  \int_{\{|h(\cdot, z)| > a \}} h^2(y, z) \, dy \geq a^2 
  \int_{\{|h(\cdot, z)| > a \}}  \, dy.
\end{eqnarray}
Property (W2) implies that $w(y,z) \le 2a e^{-\lambda_{-}z}$ for
sufficiently large $z$. Using this observation and inequality
(\ref{huliRz}), we obtain
\begin{eqnarray}
  \int_{R}^{\infty} \int_{\{|h(\cdot, z)| > a \}} e^{cz}
  w(y,z)dydz\leq 2a 
  \int_R^{\infty}\left( e^{\left(c-\lambda_{-}\right)z}
    \int_{\{|h(\cdot, z)| > a \}}  \, dy\right) dz \nonumber\\
  \leq 
  \frac{2}{ac}\int_{R}^{\infty}\left(e^{-\lambda_{-}z} \int_z^\infty
    \int_\Omega e^{c\zeta} h_{\zeta}^2 \,  dy d\zeta\right)dz
  \nonumber \\ 
  \leq 
  \frac{2}{ac}\int_{R}^{\infty}e^{-\lambda_{-}z}dz \int_R^\infty
  \int_\Omega e^{cz} h_{z}^2 \,  dy dz \nonumber \\ 
  = \frac{2}{ac\lambda_{-}} e^{-\lambda_{-} R} \int_R^\infty
  \int_\Omega e^{cz} h_{z}^2 \,  dy dz \leq 
  \frac{2}{ac\lambda_{-}} e^{-\lambda_{-} R} 
  \int_{\Sigma} e^{cz} h_{z}^2 \,  dx.
\end{eqnarray}
It then follows
\begin{eqnarray}
  \int_{\{|h| > a \}}  e^{cz} \bigl| V(w + h, y) - V'(w, y)
  h - V(w, y) - V(h, y) \bigr| dx \nonumber \\
 \leq C_4 a e^{c R} + C_5  a^{-1} e^{-{\lambda_- R 
      }} \int_\Sigma e^{cz} h_z^2 \,  dx,
\end{eqnarray}
for some $C_4, C_5 > 0$.

In view of the above, for a given $a$ one can choose $R$ sufficiently
large so that
\begin{eqnarray}
  \int_{\{|h| > a \}} e^{cz} \bigl| V(w + h, y) - V'(w, y)
  h - V(w, y) - V(h, y) \bigr| dx \nonumber \\ 
  \leq K_2 + C_6 a^\gamma \int_\Sigma e^{cz} h_z^2 \,  dx,
\end{eqnarray}
where $K_2 = K_2(a)$ and $C_6 > 0$ is a constant independent of $a$.

\paragraph{Step 4}

Finally, combining the estimates of Steps 1 and 3 and choosing $a$
small enough, we obtain
\begin{eqnarray}
  \Psi_c^w[h] \geq -K_3  + C_7 (c - c^* - \ep) \int_\Sigma e^{cz}
  h_z^2 \,  dx,
\end{eqnarray}
where $K_3, C_7 > 0$, for any given $\ep > 0$. Therefore, choosing
$\ep$ sufficiently small and also making use of the Poincar\'e
inequality \cite[Lemma 2.1]{lmn:arma08}, we arrive at
\begin{eqnarray}
  \int_\Sigma e^{cz} h_z^2 \,  dx \leq K, \qquad
  \int_\Sigma e^{cz} h^2 \, dx \leq K,
\end{eqnarray}
for some $K > 0$ on a sub-level set of $\Psi_c^w$. This, in turn,
implies that 
\begin{eqnarray}
  \int_\Sigma e^{cz} |\nabla_y h|^2 dx \leq 2 \Psi_c^w[h] +
  \int_\Sigma e^{cz} |V''(\tilde u, y)| h^2 dx \leq K',
\end{eqnarray}
for some $K' > 0$, on that sub-level set as well, and hence the
functional $\Psi_c^w$ is coercive for all $c > c^*$.
\end{proof}

\section{Existence and properties of minimizers of $\Phi_c^w$}
\label{s:exist} 

In this section we establish existence of traveling wave solutions of
(\ref{pde}) with exponential decay governed by (\ref{eq:lammin}),
which are global minimizers of the respective functionals. Let us
point out that in general this approach may not give all such
solutions of (\ref{trav}) and \eqref{bc}. Nevertheless, as we will
show later in Section \ref{sec:uniq-supercr-fronts}, these are in fact
the only such solutions for large enough values of $c$. Also, for
certain classes of nonlinearities, the same result holds for all
traveling wave solutions with $c > c_0$ (see Corollary \ref{c:kpp}).

We begin with a basic characterization of $c^*$ in \eqref{eq:cstar} in
terms of existence of minimizers of the functional $\Phi_c$ for $c >
c_0$.
\begin{Proposition}
  \label{p:cstar}
  Let $c^*$ be defined in \eqref{eq:cstar}. Then $c^* \geq c_0$, there
  is no non-trivial minimizer of $\Phi_c$ over all $u \in
  H^1_c(\Sigma)$ satisfying \eqref{bc} for any $c > c^*$, and if $c^*
  > c_0$, there exists a non-trivial minimizer of $\Phi_{c^*}$ over
  all $u \in H^1_{c^*}(\Sigma)$ satisfying \eqref{bc}.
\end{Proposition}

\begin{proof}
  Suppose first that there exists a non-trivial minimizer of $\Phi_c$
  for some $c = c^\dag > c_0$. Then it is clear that $c^* =
  c^\dag$. Indeed, by \cite[Proposition 3.2]{lmn:arma08} and the
  arguments of \cite[Proposition 6.9]{lmn:arma08} we have $\inf \Phi_c
  = 0$ for every $c \geq c^\dag$. On the other hand, by \cite[Equation
  (5.4)]{mn:cms08} we have $\inf \Phi_c < 0$ for all $0 < c < c^\dag$.

  Alternatively, if there is no non-trivial minimizer for any $c >
  c_0$, then $\inf \Phi_c = 0$ for all $c > c_0$. Indeed, if not then
  there exists a trial function $u_c \in H^1_c(\Sigma)$ with values in
  $[0, 1]$, such that $\Phi_c [u_c] < 0$, contradicting non-existence
  of minimizers guaranteed by \cite[Theorem 3.3]{mn:cms08}. On the
  other hand, arguing as in \cite[Theorem 4.2]{mn:cms08}, it is not
  difficult to see that the trial function in \cite[Equation
  (4.18)]{mn:cms08} makes $\Phi_c$ negative for any $0 < c < c_0$.
\end{proof}

The proof of existence of traveling wave solutions in the statement of
Theorem \ref{exist} follows from a sequence of propositions below.

\begin{Proposition}
  \label{p:exist}
  Let $c > c^*$ and let $w$ satisfy (W1) -- (W3). Then there exists
  $\bar u \in C^2(\Sigma) \cap C^1(\overline\Sigma) \cap
  W^{1,\infty}(\Sigma)$, with $0 < \bar u < 1$, such that
  \begin{enumerate}
  \item[(i)] $\bar u$ solves (\ref{trav}) with (\ref{bc}).

  \item[(ii)] $\bar u$ minimizes $\Phi_c^w[u]$ over all $u$ taking
    values in $[0, 1]$ and such that $u - w \in H^1_c(\Sigma)$.
  \end{enumerate}
\end{Proposition}

\begin{proof}
  By coercivity of $\Psi_c^w$ established in Proposition \ref{p:coer},
  if $(h_n)$ is a minimizing sequence from $\mathcal A_c^w$, then
  $(h_n)$ is uniformly bounded in $H^1_c(\Sigma)$. Hence there exists
  a subsequence, still labeled $(h_n)$ such that $h_n \rightharpoonup
  \bar h$ in $H^1_c(\Sigma)$. Moreover, $\bar h$ is obviously still in
  $\mathcal A_c^w$. Therefore, by lower-semicontinuity of $\Psi_c^w$
  in the weak topology of $H^1_c(\Sigma)$ proved in Proposition
  \ref{p:lsc}, the limit $\bar h$ is a minimizer of $\Psi_c^w$ in
  $\mathcal A_c^w$.  Now, since $\Phi_c^w[u] = \Psi_c^w[u - w]$, the
  function $\bar u = w + \bar h$ is also a minimizer of $\Phi_c^w$
  over all $u$, such that $0 \leq u \leq 1$ and $u - w \in
  H^1_c(\Sigma)$. This proves part (ii) of the statement.

  To prove part (i), consider the first variation of the functional
  $\Phi_c^w$ evaluated on $\bar u$ with respect to a smooth function
  $\varphi$ with compact support, vanishing at $\partial
  \Sigma_\pm$. Since $\bar u$ is a minimizer and the barriers $u = 0$
  and $u = 1$ are sub- and supersolutions of \eqref{trav}, we have
  \cite{kinderlehrer}
  \begin{eqnarray}
    0 = \delta \Phi_c^w(\bar u, \varphi) = \int_\Sigma e^{cz} \left\{
      \nabla \bar u \cdot \nabla \varphi - f(\bar u, y) \varphi
    \right\} \, dx.
  \end{eqnarray}
  By elliptic regularity theory \cite{gilbarg} $\bar u$ is a classical
  solution of (\ref{trav}), and $\bar u\in W^{1,\infty}(\Sigma)$. 
  
  Moreover, we also have $\bar u \not \equiv 0$, since $\bar u = 0$ would imply 
  $\bar h = \bar u - w = -w \in H^1_c(\Sigma)$, contradicting our
  assumptions on $w$. By strong maximum principle, we then have $0 <
  \bar u < 1$, in view of the fact that $u = 0$ and $u = 1$ are the
  sub- and super-solutions of (\ref{trav}), respectively
  \cite{gilbarg,mn:cms08}.
\end{proof}

\begin{Remark}
  We note that the proof of Proposition \ref{p:exist} does not rely on
  the precise details of the assumptions (W1)--(W3) on the function
  $w$. The main ingredient of the proof is that $w$ solves
  \eqref{trav} and \eqref{bc} and does not lie in
  $H^1_c(\Sigma)$. Therefore, with a suitable construction of $w$ our
  approach may also be applied to supercritical traveling waves in the
  case $\nu_0 = 0$, yielding existence of traveling wave solutions
  with non-exponential decay.
\end{Remark}

It is not difficult to see that the minimizer constructed in
Proposition \ref{p:exist} has the decay governed by
(\ref{eq:lammin}). We note that exponential decay of solutions of
\eqref{trav} in cylindrical domains has been extensively studied under
various assumptions on the nonlinearity, etc. (see, for instance,
\cite{berestycki92,vega93}). In particular, with minor modifications
the arguments of \cite{vega93} could be applied in our slightly more
general setting. For the sake of completeness, we present the proof
under our set of assumptions, using the techniques of \cite{mn:cms08}
and giving a slightly stronger decay estimates. We also note that our
result establishes the {\em precise} rate of exponential decay given
by \eqref{eq:lammin} for all solutions of \eqref{trav} with $c > c^*$.

\begin{Proposition}
  \label{p:decay}
  Let $c > c^*$ and let $\bar u$ be a solution of (\ref{trav}) and
  \eqref{bc}. Then
  \begin{eqnarray}
    \label{decay}
    \bar u(y, z) & = & a e^{-\lambda_-(c, \nu_0) z} \psi_0(y)(1 +
    O(e^{-\delta z})), \\
    \label{decayz}
    \bar u_z(y, z) & = & -a \lambda_-(c, \nu_0) e^{-\lambda_-(c, \nu_0)
      z} \psi_0(y)(1 +     O(e^{-\delta z})), 
  \end{eqnarray}
  uniformly in $\Omega$ for some $a \in \mathbb R$ and $\delta > 0$.
\end{Proposition}

\begin{proof}
  We reason as in the proof of \cite[Theorem 3.3(iii)]{mn:cms08},
  recalling that $c > c_0$ by Proposition \ref{p:cstar}.  By a Fourier
  expansion of \eqref{trav} in terms of $\psi_k$ and $\nu_k$ from
  \eqref{eq:psik}, after some algebra we obtain that $a_k(z) =
  \int_\Omega \bar u(y, z) \psi_k(y) \, dy \in C^2(\mathbb R)$ satisfy
  \begin{eqnarray}
    \label{eq:ak}
    a''_k + c a'_k - \nu_k a_k = g_k,
  \end{eqnarray}
  where 
  \begin{eqnarray}
    \label{eq:gk}
    g_k(z) = \int_\Omega (f_u(0, y) - f_u(\tilde u(y, z), y)) \bar u(y, z)
    \psi_k(y) dy,  
  \end{eqnarray}
  for some $0 < \tilde u < \bar u$. Using variation of parameters, for
  any $R \in \mathbb R$ we can write the solution of \eqref{eq:ak} in
  the form
  \begin{eqnarray}
    a_k(z) =  a_k^+(R) e^{-\lambda_+(c, \nu_k) (z - R)} + a_k^-(R)
    e^{-\lambda_-(c, \nu_k) (z - R)} \nonumber \\ 
    -{1 \over \sqrt{{c}^2 + 4 \nu_k}} \int_R^z  
    e^{\lambda_+(c, \nu_k) (\xi - z)} g_k(\xi) d\xi \nonumber \\ +
    {1 \over \sqrt{{c}^2 + 4 \nu_k}} \int_R^z
    e^{\lambda_-(c, \nu_k) (\xi - z)}   
    g_k(\xi) d\xi, \label{eq:contr}
  \end{eqnarray}
  where $\lambda_\pm$ are defined in \eqref{eq:lammin}. Specifically,
  since $\psi_0 > 0$ and $\bar u > 0$, for $k = 0$ we have $a_0 > 0$
  and
  \begin{eqnarray}
    \label{eq:g0}
    |g_0(z)| \leq C ||u(\cdot, z)||_{L^\infty(\Omega)}^\gamma a_0(z),
  \end{eqnarray}
  for some $C > 0$.

  Since $\bar u(\cdot, z) \to 0$ uniformly as $z \to +\infty$, by the
  definition of $a_0(z)$ and \eqref{eq:g0} we have
\begin{eqnarray}
  \label{eq:a0zinf}
  a_0(z) \to 0 \quad \text{and} \quad {g_0(z) \over a_0(z)} \to 0
  \qquad \text{as} \qquad z \to +\infty. 
\end{eqnarray}
Therefore, by \eqref{eq:contr} for every $b_0 > 0$ there exists $\bar
R_0 \in \mathbb R$, such that for every $b \in (0, b_0)$ and every
$R_0 \geq \bar R_0$ we have
\begin{eqnarray}
  \label{eq:a0est}
  a_0(R_0) \leq a_0^+(R) e^{\lambda_+(c, \nu_0) b} + a_0^-(R)
  e^{\lambda_-(c, \nu_0) b} + \frac12 \max_{R_0 \leq z \leq R} a_0(z),
\end{eqnarray}
where $R = R_0 + b$. In fact, the maximum in the last term of
\eqref{eq:a0est} is attained at $z = R_0$ due to monotonic decrease of
$a_0(z)$ for all $z \geq \bar R_0$, provided that $\bar R_0$ is large
enough. Indeed, in view of \eqref{eq:a0zinf} we can always choose
$\bar R_0$ sufficiently large, so that $\max_{z \geq \bar R_0} a_0(z)
= a_0(\bar R_0)$. Then, if $a_0(z)$ is not monotonically decreasing,
it attains a local minimum at some $z_0 \in (\bar R_0,
+\infty)$. However, by \eqref{eq:ak} with $k = 0$ and
\eqref{eq:a0zinf} we have
\begin{eqnarray}
  \label{eq:a0pp}
  a_0''(z_0) \leq \tfrac12 \nu_0 a_0(z_0) < 0,
\end{eqnarray}
contradicting minimality of $a_0(z)$ at $z_0$.  

Recalling that $b = R - R_0$, the monotonicity of $a_0(z)$ for $z \geq
\bar R_0$ and the estimate in \eqref{eq:a0est} imply that for every
$\varepsilon > 0$ there exists $b_1 > 0$ independent of $R_0$, such
that
\begin{eqnarray}
  0 & \leq & \frac12 e^{-\lambda_+(c, \nu_0) b} a_0(R_0) \leq a_0^+(R) +
  e^{-(\lambda_+(c, \nu_0) - \lambda_-(c, \nu_0)) b} a_0^-(R)
  \nonumber \\ 
  & \leq & a_0^+(R) + \varepsilon a_0^-(R) \qquad \forall R \in [R_0 +
  b_1, R_0 + 2 b_1], 
  \label{eq:apmest}
\end{eqnarray}
where the last line follows from $\lambda_+(c, \nu_0) > \lambda_-(c,
\nu_0) > 0$, choosing $b_0$ large enough. Therefore, by
\eqref{eq:contr} and, in particular, the fact that $a_0(R) = a_0^+(R)
+ a_0^-(R)$ we have
\begin{eqnarray}
  \label{eq:apm}
  \left.  {d a_0 \over dz} \right|_{z = R} & = & -\lambda_-(c, \nu_0) 
  a_0^-(R) -  \lambda_+(c, \nu_0)
  a_0^+(R)  \nonumber \\
  & \leq & -\frac12 \lambda_-(c, \nu_0) a_0(R) < 0 \qquad
  \forall R \in [R_0 + b_1, R_0 + 2 b_1],  
\end{eqnarray}
provided that $\varepsilon$ is small enough.  Now, in view of
arbitrariness of $R_0 \geq \bar R_0$, we conclude that the inequality
\eqref{eq:apm} in fact holds for all $z \geq R_1 := R_0 + b_1$, and so
$a_0(z) \leq a_0(R_1) e^{- \frac12 \lambda_-(c, \nu_0) (z - R_1)}$ for
all $z \geq R_1$. Hence, by \cite[Eq. (4.24)]{mn:cms08} we have $\bar
u(\cdot, z) \leq C e^{-\mu z}$, for some $C > 0$ and $\mu > 0$. From
now on, we can argue exactly as in the proof \cite[Theorem
3.3(iii)]{mn:cms08} to show that either $\bar u(y, z) = a
e^{-\lambda_-(c, \nu_0) z} \psi_0(y) (1 + e^{- \delta z})$ or $\bar
u(y, z) = a e^{-\lambda_+(c, \nu_0) z} \psi_0(y) (1 + e^{- \delta
  z})$, for some $a > 0$ and $\delta > 0$. Here we took into account
that the slightly stronger multiplicative decay estimate than the one
obtained in \cite{mn:cms08} follows from the $C^1$ convergence on
slices of $\Sigma$ for the difference between $\bar u$ and the leading
order term in \cite[Eq. (3.38)]{mn:cms08} and the Hopf lemma applied
to $\psi_0$.

To finish the proof of \eqref{decay}, we need to exclude the
possibility that the exponential decay rate is governed by
$\lambda_+(c, \nu_0)$. Indeed, if the latter were true, then by the
above estimates together with similar estimates for $|\nabla \bar u|$
following from \cite[Eq. (3.38)]{mn:cms08}, we find that $\bar u \in
H^1_c(\Sigma)$. Therefore, $\bar u$ satisfies hypothesis (H3) of
\cite{mn:cms08}, implying that there exists a non-trivial minimizer of
$\Phi_{c'}$ for some $c' \geq c$. Then, by the arguments in the proof
of Proposition \ref{p:cstar}, this would imply that $c^* \geq c'$,
thus leading to a contradiction.

Finally, the estimate for $\bar u_z$ in (\ref{decayz}) is obtained by
differentiating \eqref{trav} with respect to $z$ and using the same
arguments as for (\ref{decay}).
\end{proof}

Thus, according to Propositions \ref{p:exist} and \ref{p:decay} there
exist bounded solutions of (\ref{trav}) and (\ref{bc}) with the
behavior at $z = +\infty$ governed by the slowest exponentially
decaying positive solution of the respective linearized equation. Now
we will establish a number of additional properties of these
solutions, based on the minimality of $\Phi_c^w$ evaluated on $\bar
u$. In particular, we will show that these solutions have the form of
advancing fronts.

\begin{Proposition}
  Let $\bar u$ be as in Proposition \ref{p:exist}. Then
  \begin{enumerate}
  \item[(i)] $\bar u(y, z)$ is a strictly monotonically decreasing
    function of $z$. 

  \item[(ii)] $\bar u(\cdot, z) \to v$ in $C^1(\overline\Omega)$ as $z
    \to -\infty$, where $0 < v \leq 1$ is a critical point of $E$, and
    $E[v] < 0$.
  \end{enumerate}
  \label{p:mono}
\end{Proposition}

\begin{proof}
  We prove monotonicity by using a monotone decreasing rearrangement
  of $\bar u(y, \cdot)$, following an idea of Heinze
  \cite{heinze}. Introducing a new variable $\zeta = e^{cz}$, we can
  write the functional $\Phi_c^w$ as
  \begin{eqnarray} \label{eq:2} %
    \Phi_c^w[\bar u] = {1 \over c} \int_0^\infty \int_\Omega \Biggl(
    {c^2 \zeta^2 \bar u_\zeta^2 \over 2} + {|\nabla_y \bar u|^2 \over
      2} + V(\bar u, y) \hspace{1.8cm} \nonumber \\
    - {c^2 \zeta^2 w_\zeta^2 \over 2} - {|\nabla_y w|^2 \over 2} -
    V(w, y) \Biggr) dy \, d\zeta.
  \end{eqnarray}
  Let us now perform a one-dimensional monotone decreasing
  rearrangement of $\bar u$ in $\zeta$ for $y$ fixed (see
  \cite[Section II.3]{kawohl} for details). If $\bar u^*: \Sigma \to
  [0, 1]$ is the result of this rearrangement, then $\bar u^*(\cdot,
  z) = \bar u(\cdot, z)$ for all $z \geq R$, with $R$ sufficiently
  large, in view of monotonicity of $\bar u$ guaranteed by
  (\ref{decayz}) in Proposition \ref{p:decay}. As a consequence, the
  function $\bar u^*$ belongs to the admissible class for
  $\Phi_c^w$. On the other hand, by \cite[Proposition 12]{bertsch07}
  and \cite[Lemma 2.6 and Remark 2.34]{kawohl} we know $\Phi_c[\bar
  u^*] < \Phi_c[\bar u]$, unless $\bar u = \bar u^*$, contradicting
  the minimizing property of $\bar u$. Then, by strong maximum
  principle the minimizer $\bar u(y, z)$ is, in fact, a strictly
  decreasing function of $z$. This proves part (i) of the statement.

  Once the monotonic decrease of $\bar u$ is established, by
  boundedness of $\bar u$ there exists a limit of $\bar u(\cdot, z)$
  as $z \to -\infty$. The remainder of part (ii) of the statement then
  follows exactly as in \cite[Proposition 6.6 and Corollary
  6.7]{lmn:arma08}.
\end{proof}

\begin{proof}[Proof of Theorem 1]
  The result of parts (i)--(iii) of Theorem \ref{exist} follow by
  combining the results of Propositions \ref{p:wexists},
  \ref{p:exist}, \ref{p:decay} and \ref{p:mono} and taking into
  account translational invariance in $z$, which allows us to consider
  only small values of $a > 0$. Also, in (ii) the arguments in the
  proof of \cite[Theorem 3.3]{mn:cms08} are used to establish the
  estimate on $|\nabla \bar u|$. The result of part (iv) follows
  exactly as in \cite[Theorem 3.3]{mn:cms08}. To prove part (v), note
  that if $u$ satisfies \eqref{u01} and $u - \tilde w \in
  H^1_c(\Sigma)$, then by the assumption on $\tilde w$ we have that
  $\Phi_c^{\tilde w}[u]$ is well-defined, and by minimality of
  $\Phi_c^w[\bar u]$ we have
  \begin{eqnarray}
    \label{eq:www}
    \Phi_c^{\tilde w}[\bar u] = \Phi_c^w[\bar u] + \Phi_c^{\tilde w}[w] \leq
    \Phi_c^w[u] + \Phi_c^{\tilde w}[w] \leq \Phi_c^{\tilde w}[u],
  \end{eqnarray}
  which proves the claim.
\end{proof}


\section{Global uniqueness of supercritical fronts}
\label{sec:uniq-supercr-fronts}

Let us now discuss the question whether the obtained traveling wave
solutions are unique among all traveling waves with speed $c$, up to
translations.

First of all, by Proposition \ref{p:decay} all bounded positive
solutions of (\ref{trav}) and \eqref{bc} with $c>c^*$ must have the
asymptotic decay specified by (\ref{eq:lammin}), which is the same
decay of the minimizer $\bar u$ of $\Phi_c^w$.
Using the sliding domains method of Berestycki and Nirenberg
\cite{berestycki91}, it is possible to show that any two traveling
front solutions sandwiched between two limiting equilibria must
coincide up to translation (see \cite[Theorem
5.1]{vega93}). Therefore, in order to have non-uniqueness, two
different solutions with, say, $a = 1$ fixed in (\ref{eq:w0}), should
have different limiting behavior at $z = -\infty$.

\begin{figure}
\hspace{0.2cm} {\bf a)} \hspace{5.6cm} {\bf b)}
\vspace{-8mm}
\begin{center}
\includegraphics[width=5cm]{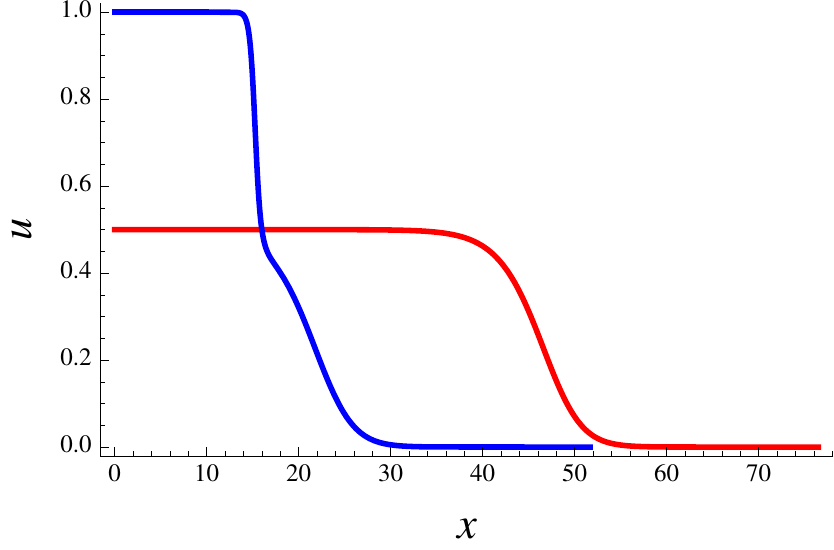}
\hspace{1cm}
\includegraphics[width=5cm]{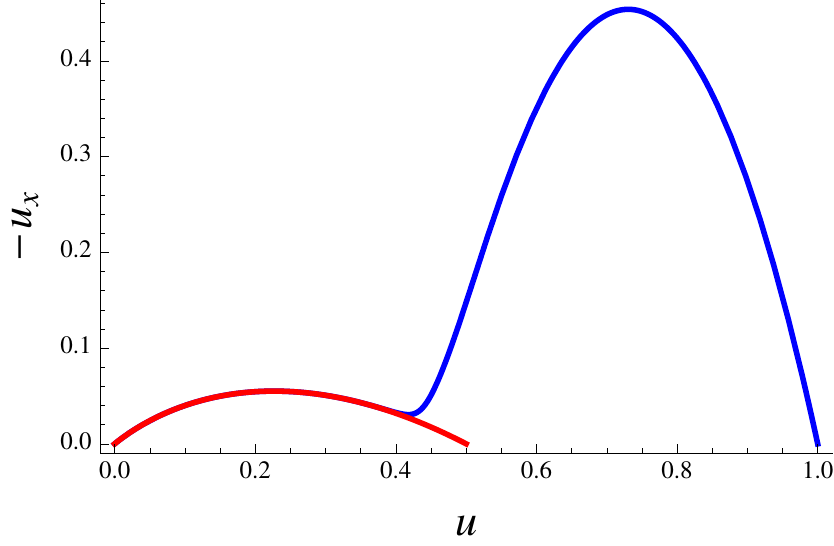}
\end{center}
\caption{A numerical example of non-uniqueness of fronts (see text for
  details).}
  \label{fig:nonu}
\end{figure}

\paragraph{Example of non-uniqueness. } Clearly, uniqueness of the
traveling front solutions with prescribed exponential decay cannot
hold true in general, even in dimension one.  Let us demonstrate such
non-uniqueness by an explicit example. For simplicity, we present a
numerical demonstration of existence of two traveling wave solutions
with the same value of the velocity $c$, the result is then made
rigorous by applying the Theorems of Sec. \ref{s:mr}. Namely, let us
consider the problem (\ref{trav}) in the one-dimensional setting with
the nonlinear term $f$ given by
\begin{eqnarray}\label{eq:fe}
  f(u)=\left\{ \begin{array}{ll}
      u(1-2u), & 0 \le u< \tfrac12,  \\
      (2u-1)(1-u)(40u-21), & \tfrac12 \le u \le 1.
\end{array} \right.
\end{eqnarray}
As can be easily verified, this function satisfies conditions (H1),
(H2) and (U).  One may expect that for a range of velocities $c$ close
to $c_0 = 2$ the problem (\ref{trav}) with nonlinearity (\ref{eq:fe})
has two traveling front solutions: one connecting $u = 0$ and $u =
\tfrac12$, and another connecting $u = 0$ with $u = 1$. These two
solutions for $c = 2.25$ are presented in Fig. \ref{fig:nonu}, where
Fig. \ref{fig:nonu}(a) shows the traveling wave profiles and
Fig. \ref{fig:nonu}(b) shows the corresponding phase portraits. From
the exact traveling wave solution $u(z) = \frac34 - \frac14 \tanh
\left( \sqrt{5 \over 2} \, z \right)$ connecting $v_1 = \frac12$ with
$v = 1$, we find that $c^\dag_{v_1} = {9 \over \sqrt{10}} \approx
2.846$ (see Sec. \ref{s:mr} for the details of the notations). Also,
by \cite[Theorem 5.1]{lmn:cpam04} we have $c^* = c^*_1 =
2$. Therefore, since $c^\dag_{v_1} > c^*$, there are two traveling
wave solutions with speed $c = 2.25$ by Remark \ref{r:2}, one
connecting zero to $v_1$ and the other connecting zero to $v$. This is
what we see numerically in Fig. \ref{fig:nonu}. Note that the value of
$c^\sharp$ from \eqref{eq:c1} is easily found to be $c^\sharp =
\sqrt{127 \over 5} \approx 5.04$.

Let us now present some heuristic arguments as to why the solutions of
\eqref{trav} and \eqref{bc} are expected to be unique for large enough
$c$. Suppose that the nonlinearity $f$ becomes {\em strictly} linear
for sufficiently small $u$, i.e. we have $f(u, y) \equiv f_u(0, y) u$
for all $u \leq \delta$, with some $\delta > 0$. In this case,
clearly, the decay of any positive solution of (\ref{eq:w}) must be
precisely
\begin{eqnarray}
  \label{eq:declin}
  w(y, z) = a e^{-\lambda_-(c, \nu_0) z} \, \psi_0(y) + \sum_{k =
    0}^\infty a_k^+ e^{-\lambda_+(c, \nu_k)} \, \psi_k(y), \qquad z
  \geq z_0,
\end{eqnarray}
for some $z_0 \gg 1$. Here we took into account that when $\nu_0 <
\nu_k \leq 0$, if ever, we have $\lambda_-(c, \nu_k) < \lambda_-(c,
\nu_0)$, hence the presence of a term proportional to
$e^{-\lambda_-(c, \nu_k) z} \psi_k(y)$ would contradict positivity of
the solution.  Now, fixing the value of $a$, one can see that any two
solutions of the form (\ref{eq:declin}) differ by a function that lies
in $H^1_c(\Sigma)$, since all the terms in the series in
(\ref{eq:declin}) except the first one decay exponentially with the
rate strictly greater than $e^{-cz/2}$. On the other hand, it is easy
to see that for $c > c^\sharp$, where $c^\sharp$ is defined in
(\ref{eq:c1}), the functional $\Psi_c^w$ is strictly convex in
$H^1_c(\Sigma)$. Therefore, for every value of $a$ and every $c >
c^\sharp$ there exists a {\em unique} positive solution of
(\ref{trav}). This seems rather surprising, since, on one hand,
uniqueness in this case is a consequence of convexity of $\Psi_c^w$
for $c > c^\sharp$, which, in some sense, is dictated by what happens
in the ``core'' of the traveling front and thus represents a nonlinear
feature of the problem. On the other hand, uniqueness also relies on
the fact that, regardless of the value of $c$, any two positive
solutions with the same asymptotic decay have difference that lies in
$H^1_c(\Sigma)$, which is a property of the behavior far ahead of the
front and represents a linear feature of the problem. From the
variational point of view, uniqueness of solution is also rather
surprising, since, in view of Remark \ref{r:walot}, there are
infinitely many auxiliary functions $w$ that can be used in the
construction the functional $\Phi_c^w$ and hence, potentially,
different minimizers, depending on the choice of $w$.
 
These heuristic arguments turn out to work in the case of generic
nonlinearities in \eqref{trav}, with the result stated in Theorem
\ref{unique}. The proof of this theorem follows from a sequence of
lemmas below.

\begin{Lemma}
  \label{l:c1}
  We have $c^\sharp \geq c^*$.
\end{Lemma}

\begin{proof}
  First, by the definition of $c^\sharp$ we trivially have $c^\sharp
  \geq c_0$. Now, take any $u \in H^1_{c^\sharp}(\Sigma)$. Then by
  hypothesis (H2) there exists $\tilde u : \Sigma \to \mathbb R$, with
  $0 < \tilde u < u$, such that
  \begin{eqnarray}
    \label{eq:16}
    \Phi_{c^\sharp}[u] = \int_\Sigma e^{c^\sharp z} \left( \frac12 |\nabla u|^2 -
      f_u(\tilde u, y) u^2 \right) dx. 
  \end{eqnarray}
  Applying the definitions of $q$ and $\hat \nu_0$ in \eqref{eq:hnu0},
  and the Poincar\'e inequality \cite[Lemma 2.1]{lmn:arma08}, we
  obtain
  \begin{eqnarray}
    \label{eq:21}
    \Phi_{c^\sharp}[u] \geq \frac12 \int_\Sigma e^{c^\sharp z} \left( u_z^2 +
      |\nabla_y u|^2 - q(y) u^2 \right) \, dx \nonumber \\  
    \geq \frac12 \int_\Sigma e^{c^\sharp z} \left( {{c^\sharp}^2 \over 4} +  
      \hat \nu_0 \right) u^2 \, dx \geq 0,
  \end{eqnarray}
  implying that by \eqref{eq:cstar} we have $c^* \leq c^\sharp$.
\end{proof}

\begin{Lemma}
  \label{l:decay}
  Let $c > c_0$ and let $\bar u_1$ and $\bar u_2$ be two solutions of
  (\ref{trav}) and \eqref{bc}, and whose decay at $z = +\infty$ is
  given by \eqref{decay}. Then $\bar u_2 - \bar u_1 \in
  H^1_c(\Sigma)$.
\end{Lemma}

\begin{proof}
  Introduce $h = \bar u_2 - \bar u_1$, which satisfies
  \begin{eqnarray}
    \label{eq:3}
    \Delta h + c h_z + f_u(0, y) h = ( f_u(0, y) - f_u(\tilde u, y))
    h,  
  \end{eqnarray}
  for some $\tilde u$ bounded between $\bar u_1$ and $\bar u_2$. In
  view of Proposition \ref{p:decay}, we have
  \begin{eqnarray}
    \label{dectil}
    0 \leq \tilde u(y, z) \leq \min \{1, C e^{-\lambda_-(c, \nu_0) z}  \}, 
  \end{eqnarray}
  for some $C > 0$ and all $(y, z) \in \Sigma$.

  Now, let $\eta \in C^1(\mathbb R)$ be a cutoff function with the
  following properties: $\eta(z)$ is a non-increasing function of $z$
  with super-exponential decay at $+\infty$, and $\eta(z) = 1$ for $z
  \leq 0$. For fixed $\ep > 0$ and $c' > 0$ let us multiply
  \eqref{eq:3} with $ e^{c' z} h \eta_\ep^2$, where $\eta_\ep(z) =
  \eta(\ep z)$, and integrate over $\Sigma$. After a few integrations
  by parts, we obtain
  \begin{eqnarray}
    \label{eq:1}
    \int_\Sigma e^{c' z} \left( |\nabla h|^2 - {c' (c' - c) \over 2}
      h^2 - f_u(0, y) h^2 \right)  \eta_\ep^2 \, dx  \nonumber \\ 
    +  \int_\Sigma e^{c' z} \left( 2 h h_z \eta_\ep \eta'_\ep -
      (c' - c) h^2 \eta_\ep \eta'_\ep \right) \, dx \nonumber \\ 
    =  \int_\Sigma e^{c' z}  (f_u(\tilde u,
    y) - f_u(0, y)) h^2 \eta_\ep^2 \, dx,
  \end{eqnarray}
  where $\eta_\ep' = d \eta_\ep / dz$.
  
  We are now going to estimate the different terms appearing in
  \eqref{eq:1}. Let us begin by noting that by hypothesis (H2) and
  \eqref{dectil} we have $|f_u(\tilde u, y) - f_u(0, y)| \leq C
  |\tilde u(y, z)|^\gamma \leq C' (1 + e^{\gamma \lambda_-(c, \nu_0)
    z})^{-1}$ for some $C, C' > 0$ and all $(y, z) \in
  \Sigma$. Therefore,
  \begin{eqnarray}
    \label{eq:4}
     \int_\Sigma e^{c' z}  (f_u(\tilde u, y) - f_u(0, y)) h^2
     \eta_\ep^2 
     \, dx \leq C' \int_\Sigma e^{c' z} (1 + e^{\gamma \lambda_-(c,
       \nu_0) z} )^{-1} h^2 \eta_\ep^2 \,  dx. 
  \end{eqnarray}
  On the other hand, since $h \eta_\ep \in H^1_{c'}(\Sigma)$, by
  Poincar\'e inequality \cite[Lemma 2.1]{lmn:arma08} it holds
  \begin{eqnarray}
    \label{eq:5}
    \int_\Sigma e^{c' z} h_z^2 \eta_\ep^2 \, dx \geq {{c'}^2
      \over 4} \int_\Sigma e^{c' z} h^2 \eta_\ep^2 \, dx \hspace{3cm} 
    \nonumber \\ 
    - 2 \int_\Sigma e^{c' z} h h_z \eta_\ep \eta'_\ep \, dx - 
    \int_\Sigma e^{c' z}h^2  |\eta'_\ep|^2 \, dx.
  \end{eqnarray}
  Combining this with the previous estimate and using the variational
  characterization of the lowest eigenvalue $\nu_0$ in \eqref{eq:nu0},
  from \eqref{eq:1} we obtain
  \begin{eqnarray}
    \label{eq:6}
    \int_\Sigma e^{c' z} \left( {c c' \over 2} - {{c'}^2 \over 4} +
      \nu_0 \right) h^2 \eta_\ep^2 \, dx \hspace{4cm} \nonumber \\ 
    -  |c' - c|  \int_\Sigma
    e^{c' z}  h^2 \eta_\ep |\eta_\ep'| \, dx - \int_\Sigma e^{c'
      z} h^2  |\eta_\ep'|^2 \, dx \nonumber \\ 
    \leq  C \int_\Sigma e^{c' z} (1 + e^{\gamma \lambda_-(c, \nu_0) z}
    )^{-1} h^2 \eta_\ep^2 \,  dx,
  \end{eqnarray}
  for some $C > 0$ independent of $\ep$.
  
  To proceed, we need to choose the cutoff function $\eta$ in such a
  way that the terms involving $\eta_\ep'$ can be bounded by the terms
  involving $\eta_\ep$ uniformly in $\ep$ as $\ep \to 0$. This can be
  achieved by setting $\eta(z) = \exp(-\tfrac12 z^2)$ for $z >
  0$. Indeed, consider the last term in the left-hand side of
  (\ref{eq:6}). We can write
  \begin{eqnarray}
    \label{eq:7}
    \int_\Sigma e^{c' z} h^2  |\eta_\ep'|^2 \, dx \leq \alpha^2
    \int_\Sigma e^{c' z}  h^2 \eta_\ep^2 \, dx +  \int_{ \{|\eta'_\ep|
      \geq \alpha \eta \}} e^{c' z} h^2 |\eta_\ep'|^2 \, dx,
  \end{eqnarray}
  for any $\alpha > 0$. From our choice of $\eta$, the second integral
  is over all $z \geq \alpha \ep^{-2}$. Therefore, assuming $h \in
  L^2_{c''}(\Sigma)$ for some $c'' \in (0, c')$, we have
  \begin{eqnarray}
    \label{eq:8}
    \int_\Sigma e^{c' z} h^2  |\eta_\ep'|^2 \, dx \leq \alpha^2
    \int_\Sigma e^{c' z}  h^2 \eta_\ep^2 \, dx +  \int_{\alpha
      \ep^{-2}}^{\infty} \int_\Omega e^{c' z} h^2 |\eta_\ep'|^2 \, dy
    dz \nonumber \\ 
    \leq \alpha^2 \int_\Sigma e^{c' z}  h^2 \eta_\ep^2 \, dx + C_\ep
    \int_\Sigma  e^{c'' z}  h^2 \, dx,
  \end{eqnarray}
  where
  \begin{eqnarray}
    \label{eq:11}
    C_\ep = \sup_{z \geq \alpha \ep^{-2}} \left( e^{(c' - c'') z}
      |\eta_\ep'(z)|^2 \right).
  \end{eqnarray}
  Similarly, 
  \begin{eqnarray}
    \label{eq:10}
    \int_\Sigma e^{c' z}  h^2 \eta_\ep |\eta_\ep'| \, dx \leq \alpha
    \int_\Sigma e^{c' z}  h^2 \eta_\ep^2 \, dx + \alpha^{-1} C_\ep
    \int_\Sigma e^{c'' z}  h^2 \, dx. 
  \end{eqnarray}
  Therefore, we can rewrite (\ref{eq:6}) in the following form (also
  taking (\ref{eq:c0}) into consideration)
  \begin{eqnarray}
    \label{eq:12}
    \int_\Sigma e^{c' z} \left( {c c' \over 2} - {{c'}^2 \over 4} -
      {c_0^2 \over 4} - \alpha |c' - c| - \alpha^2 \right) h^2
    \eta_\ep^2 \, dx \nonumber \\  
    \leq  C \sup_{z \in \mathbb R} \left( {e^{(c' - c'') z}
        \over 1 + e^{\gamma \lambda_-(c, \nu_0) z} } \right)
    \int_\Sigma e^{c'' z}  h^2 \,  dx \nonumber \\ 
    + C_\ep \left( 1 +  {|c' - c| \over \alpha} \right) \int_\Sigma
    e^{c'' z} h^2  \, dx , 
  \end{eqnarray}
  for some $C > 0$ independent of $\ep$.

  Now, by an explicit computation, we have
  \begin{eqnarray}
    \label{eq:9}
    C_\ep = \sup_{\zeta \geq 1} \left( \alpha^2
      \zeta^2 \exp ( -\ep^{-2} (\alpha^2 \zeta^2  - (c' - c'') \alpha
      \zeta) )  \right) \leq \alpha^2,
  \end{eqnarray}
  as soon as $c' - c'' < \alpha$ and $\ep$ is small
  enough. Introducing the constants
  \begin{eqnarray}
    \label{eq:13}
    c_\pm = c \pm \sqrt{c^2 - c_0^2},
  \end{eqnarray}
  and using \eqref{eq:9}, we then find that \eqref{eq:12} is
  equivalent to
  \begin{eqnarray}
    \label{eq:14}
    \Bigl\{ (c' - c_-) (c_+ - c') - 8 \alpha |c - c'| \Bigr\} 
    \int_\Sigma e^{c'  z} h^2 \eta_\ep^2 \, dx
    \leq  C  \int_\Sigma e^{c'' z} h^2 \, dx,
  \end{eqnarray}
  with $C > 0$ independent of $\ep$, whenever $c' - c'' < \alpha \leq
  \min( \gamma \lambda_-(c, \nu_0), |c' - c|)$ and $\ep$ is
  sufficiently small.

  We now note that by Proposition \ref{p:decay} we have $h \in
  L^2_{c''}(\Sigma)$ for some $c'' \in (c_-, c_+)$. Therefore, we can
  choose $c' = c'' + \tfrac12 \alpha$, with $\alpha$ so small that the
  quantity in the curly brackets in (\ref{eq:14}) is positive. Then,
  passing to the limit $\ep \to 0$ in (\ref{eq:14}), by monotone
  convergence theorem we can conclude that, in fact, we have $h \in
  L^2_{c'}(\Sigma)$. Applying the same argument to (\ref{eq:1}) we
  find that also $h \in H^1_{c'}(\Sigma)$, with $c' \in (c'',
  c_+)$. Finally, iterating this argument, we find that
  \begin{eqnarray}
    \label{eq:15}
    h \in H^1_{c'}(\Sigma), \qquad \forall c' < c_+.
  \end{eqnarray}
  In particular, $h \in H^1_c(\Sigma)$, which completes the proof.
\end{proof}

\begin{Lemma}
  \label{l:convex}
  Let $c > c^\sharp$. Then the functional $\Psi_c^w$ is strictly
  convex on $\mathcal A_c^w$.
\end{Lemma}
\begin{proof}
  First of all, let us point out that due to the presence of the
  exponential weight the functional $\Psi_c^w$ is not twice
  continuously differentiable in $H^1_c(\Sigma)$, so the convexity of
  the functional cannot be concluded by evaluating the second
  variation of $\Psi_c^w$. Nevertheless, what we will show below is
  that the functional $\Psi_c^w$ is {\em strongly} convex
  \cite{rockafellar}, i.e., there exists $\sigma > 0$ such that for
  every $h_1, h_2 \in \mathcal A_c^w$ it holds
  \begin{eqnarray}
    \label{eq:17}
    \Psi_c^w[h_2] \geq \Psi_c^w[h_1] + \delta \Psi_c^w(h_1, h_2 - h_1)
    + {\sigma \over 2} ||h_2 - h_1||^2_{L^2_c(\Sigma)}. 
  \end{eqnarray}
  Indeed, by hypothesis (H2) we have
  \begin{eqnarray}
    \label{eq:18}
    \Psi_c^w[h_2] & = & \Psi_c^w[h_1] + \delta \Psi_c^w(h_1, h_2 -
    h_1) +  \frac12 \int_\Sigma e^{cz} |\partial_z (h_2 - h_1)|^2 dx 
    \nonumber \\ 
    & + & \frac12 \int_\Sigma e^{cz} \left( |\nabla_y(h_2 - h_1)|^2 -
      f_u(\tilde u, y) (h_2 - h_1)^2  \right) dx,
  \end{eqnarray}
  for some $0 \leq \tilde u \leq 1$. Now, by Poincar\'e inequality
  \cite[Lemma 2.1]{lmn:arma08}, we have
  \begin{eqnarray}
    \label{eq:19}
    \frac12 \int_\Sigma e^{cz} |\partial_z (h_2 - h_1)|^2 dx  \geq {c^2 \over
      8} ||h_2 - h_1||^2_{L^2_c(\Sigma)},
  \end{eqnarray}
  and by the definition of $\hat \nu_0$ in \eqref{eq:hnu0}
  \begin{eqnarray}
    \label{eq:20}
    \frac12 \int_\Sigma e^{cz} \left( |\nabla_y(h_2 - h_1)|^2 -
      f_u(\tilde u, y) (h_2 - h_1)^2  \right) dx \geq {\hat \nu_0
      \over 2} ||h_2 - h_1||^2_{L^2_c(\Sigma)}. 
  \end{eqnarray}
  Strong convexity then follows with $\sigma = \tfrac14 c^2 + \hat
  \nu_0 > 0$, in view of $c > c^\sharp$ and implies strict convexity
  \cite[Chapter 12.H$^*$]{rockafellar}.
\end{proof}

\begin{proof}[Proof of Theorem 3] Let $\bar u_1$ and $\bar u_2$ be as
  in Lemma \ref{l:decay}. Without loss of generality we may assume
  that $a$ is sufficiently small in \eqref{decay}. Let us define $w$
  to be a suitably truncated version of $\bar u_1$. Then, by
  hypothesis (H2) the functional $\Psi_c^w$ is of class $C^1$ and by
  Lemma \ref{l:convex} is also strictly convex on $\mathcal
  A_c^w$. Then $h_1 := \bar u_1 - w$ belongs to $\mathcal A_c^w$. In
  fact, the function $h_2 := \bar u_2 - w$ also belongs to $\mathcal
  A_c^w$, in view of Lemma \ref{l:decay}. On the other hand, since the
  equivalent functional $\Phi_c^w$ is also strictly convex, it has at
  most one critical point. Therefore, since $\bar u_1$ and $\bar u_2$
  solve the Euler-Lagrange equation \eqref{trav} for $\Phi_c^w$ and
  are, therefore, critical points of $\Phi_c^w$, we conclude that
  $\bar u_1 = \bar u_2$.
\end{proof}

\section{Global asymptotic stability of fast fronts}
\label{sec:glob-asympt-stab}

We now investigate the stability properties of traveling wave
solutions with speed $c > c^*$ obtained in the preceding
sections. Let us point out in the first place that one cannot expect
these traveling wave solutions to be stable in the usual sense (say,
with respect to small $L^\infty$ perturbations). Indeed, under
hypotheses (H1), (H2) and (U) a truncation of the traveling wave
solution with any $c > c^*$ by a cutoff function that vanishes
identically for all $z \geq R \gg 1$ would not converge to the
corresponding traveling wave, but will instead propagate with the
asymptotic speed $c^*$ of the critical front \cite[Theorem
5.8]{mn:cms08}. Similarly, modifying the rate of exponential decay in
front of the traveling wave solution with speed $c > c^*$ may lead to
acceleration or slowdown of the front (see
e.g. \cite{sattinger76,bramson88,mallordy95}), or even lead to
irregular behavior \cite{yanagida07}. Therefore, stability of these
traveling wave solutions should be considered with respect to suitable
classes of perturbations \cite{sattinger76,sattinger77}. Here we show
global stability of the (unique) supercritical traveling waves with $c
> c^\sharp$ under perturbations lying in the weighted space
$L^2_{c'}(\Sigma)$ with suitably chosen $c' > 0$.

We begin with a general existence result for the solutions of
\eqref{pde} of the form
\begin{eqnarray}
  \label{eq:uw}
  u(y, z, t) = \bar u(y, z - c t) + w(y, z - ct, t),
\end{eqnarray}
where $\bar u$ is given by Theorem \ref{exist} with some fixed $c >
c^*$, and $w$ lying in an appropriate exponentially weighted
space. The function $w$ satisfies the following parabolic equation
\begin{eqnarray}
  \label{eq:wt}
  w_t = \Delta w + c w_z + f(\bar u + w, y) - f(\bar u, y).
\end{eqnarray}
Easily adapting \cite[Proposition 3.1]{mn:sima11} (see also
\cite{mn:cms08}), we can obtain the following basic existence and
regularity result for the initial value problem associated with
\eqref{eq:wt}:

\begin{Proposition}
  \label{p:wexist}
  For any $c' > 0$, let $w_0 \in L^2_{c'}(\Sigma)$ and let $0 \leq
  \bar u + w \leq 1$. Then there exists a unique solution $w \in
  C^\alpha((0, \infty); H^2_{c'}(\Sigma)) \cap C^{1,\alpha}((0,
  +\infty); L^2_{c'}(\Sigma)) \cap C([0, +\infty); L^2_{c'}(\Sigma))$,
  for all $\alpha \in (0, 1)$, of \eqref{eq:wt} and the same boundary
  conditions as in \eqref{bc}, such that $w(x, 0) = w_0(x)$.
\end{Proposition}

We are now in a position to complete the proof of Theorem
\ref{converge}, following the method introduced in \cite{m:dcdsb04}.

\begin{proof}[Proof of Theorem 3]
  Without loss of generality, we may assume that $c' < 2 \lambda_+(c,
  \hat\nu_0)$, where $\lambda_\pm(c, \nu)$ are defined in
  \eqref{eq:lammin}. The result in Proposition \ref{p:wexist} implies
  that we can multiply \eqref{eq:wt} by $e^{c'z} w$ and integrate over
  $\Sigma$. Using hypothesis (H2) and the definition of $q(y)$ in
  \eqref{eq:hnu0}, we have
  \begin{eqnarray}
    \label{eq:24}
    \int_\Sigma e^{c'z} w w_t \, dx \leq  \int_\Sigma e^{c'z} \left(
      \Delta w + c w_z - q(y) w \right) w \, dx.
  \end{eqnarray}
  Then, after a number of integrations by parts, we obtain
  \begin{eqnarray}
    \label{eq:23}
    {d \over dt} \int_\Sigma e^{c'z} w^2 \, dx \leq - 2 \int_\Sigma
    e^{c'z} \left( w_z^2 + {c' (c - c') \over 2} w^2 + |\nabla_y w|^2
      - q(y) w^2 \right) dx.
  \end{eqnarray}
  Applying the Poincar\'e inequality \cite[Lemma 2.1]{lmn:arma08} and
  the definition of $\hat \nu_0$ in \eqref{eq:hnu0}, we then arrive at
  \begin{eqnarray}
    \label{eq:25}
    {d \over dt} \int_\Sigma e^{c'z} w^2 dx \leq 2 \left(
      {{c'}^2 \over 4} - {c c' \over 2} - \hat \nu_0 \right)
    \int_\Sigma e^{c'z} w^2 \, dx \nonumber \\
    =  \frac12 \{c' - 2 \lambda_-(c, \hat\nu_0)\} \{c' - 2
    \lambda_+(c, \hat\nu_0)\}  \int_\Sigma e^{c'z} w^2 \, dx. 
  \end{eqnarray}
  Therefore, since $c' \in (2 \lambda_-(c, \hat\nu_0), 2 \lambda_+(c,
  \hat\nu_0))$, the coefficient in front of the integral in the
  right-hand side of \eqref{eq:25} is negative, which implies
  exponential decay of the $L^2_{c'}$-norm of $w$.
\end{proof}

\section*{Acknowledgments}
\label{sec:acknowledgements}

The work of PVG was supported, in part, by the United States--Israel
Binational Science Foundation grant 2006-151. CBM acknowledges partial
support by NSF via grants DMS-0718027 and DMS-0908279.

\bibliographystyle{unsrt}
 
\bibliography{../mura,../nonlin,../bio,../stat}

\end{document}